\documentclass{article}

\newcommand{\R}{I\!\! R}

\newcommand{\Z}{I\!\! Z}
\newcommand{\N}{I\!\! N}

\begin{document}

\centerline {\bf Composition series of affine manifolds and
$n-$gerbes.}

\bigskip

\bigskip

\centerline {by}

\medskip

\centerline{\bf Tsemo Aristide.}

\centerline { The Abdus Salam center for theoretical Physics}

\centerline { Strada Costiera 11}

\centerline { 34014 Trieste, Italy.}

\centerline {tsemo@ictp.trieste.it}

\bigskip

\centerline{\bf Abstract.}

\medskip

{\it In this paper, we study $n-$composition series of affine
manifolds these are  sequences $(M_n,\nabla_{M_n})\rightarrow
(M_{n-1},\nabla_{M_{n-1}})\rightarrow...(M_1,\nabla_{M_1})$, where
each affine map $f_i:(M_{i+1},\nabla_{M_{i+1}})\rightarrow
(M_i,\nabla_{M_i})$ is  surjective.
 One composition series are classified by gerbe theory. It is
 natural to think that $n-$composition series must be classified
 by $n-$gerbe theory. In the last section of this paper we propose
 a notion of abelian $n-$gerbe theory.}

\bigskip

\bigskip

{\bf Introduction.}

\bigskip

An affine bundle is a surjective affine map between affine
manifolds. A composition serie of affine manifolds is a sequence
$(M_n,\nabla_n)\rightarrow
(M_{n-1},\nabla_{M_{n-1}})\rightarrow...\rightarrow
(M_2,\nabla_{M_2})\rightarrow (M_1,\nabla_{M_1})$, where each map
$f_i:(M_{i+1},\nabla_{M_{i+1}})\rightarrow (M_i,\nabla_{M_i})$ is
an affine bundle.

 When the source space $M$, of an affine bundle is
compact, it becomes
 a locally
trivial differentiable bundle by a well-known Ehresmann theorem ( see [God]
theorem 2.11 p.16).
Let $B$ and $F$ be respectively the base and the fiber spaces of an affine
bundle with compact total space.
If moreover the second homotopy group of $B$ is trivial, then by the
 Serre bundle theorem, one deduces that
the first homotopy group $\pi_1(M)$ of $M$ is an extension of $\pi_1(B)$
 by $\pi_1(F)$.
In particular this happens when $(B,\nabla_B)$ is geodesically complete.

Auslander has conjectured that the fundamental group of a compact
geodesically complete
 affine manifold is polycyclic. The existence of a non trivial
affine map, on a finite cyclic galoisian  cover of a $n-$compact
and complete affine manifold $(n>2)$ endowed with a complete
structure eventually different from the pull-back, implies the
Auslander conjecture [T4]. The classification of affine bundles
whose total spaces are compact and complete and more generally of
composition series of affine manifolds
 will conjecturally
allow us to know all compact and complete affine manifolds, up to a finite cover,
as we know the $2-$closed and complete affine manifolds.

 The main goal of this paper is to study composition series of affine manifolds.
First we study affine bundles.

Let $\pi_1(F)$ and $\pi_1(B)$ be two groups. Write ${\R}^{m+l}={\R}^m\oplus
{\R}^l$. We denote  by $Aff({\R}^m,{\R}^l)$ the group of affine maps
of ${\R}^{m+l}$ which preserve ${\R}^l$, and by
$Aff_I({\R}^m,{\R}^l)$ the subgroup of $Aff({\R}^m,{\R}^l)$
whose restriction on ${\R}^m$ is the identity.

An algebra problem related to this classification problem of
affine bundles is the following:

Given two representations $\pi_1(F)\rightarrow Aff({\R}^l)$, and of
$\pi_1(B)\rightarrow Aff({\R}^m)$, classify all  commutative diagrams:

$$ \matrix {& 1  \rightarrow &
\pi_1(F) &\rightarrow & \pi_1(M)  \rightarrow&\pi_1(B)&
\rightarrow  1 \cr& \ \  & \downarrow &\ \ & \downarrow &\ \
\downarrow \cr & 1\rightarrow &Aff_I({\R}^m,{\R}^l) &\rightarrow
&Aff({\R}^{m},{\R}^l)& \rightarrow Aff({\R}^m)& \rightarrow 1},$$
 where the first line is an exact sequence.

There are many ways to solve the classification   problem of
affine bundles. First, we give a classification, of affinely
locally trivial affine bundles, (see definition 2.2) after we
solve the general case.

Let's present the classification of affinely locally trivial
affine bundles.
Let us consider an affinely locally trivial affine bundle with
base $(B,\nabla_B)$ and
typical
fiber $(F,\nabla_F)$. We denote by $T_F$ the group of  translations of
$(F,\nabla_F)$ (see section 3).
The affine bundle $f$ gives rise to  representations
$\pi_f:\pi_1(B)\rightarrow Aff(F,\nabla_F)/T_F$,
 $\pi_f':\pi_1(B)\rightarrow Gl(T_F)$,
and to a flat bundle $\hat T_F$, with typical fiber $T_F$ over $B$
associated to $\pi_f'$.
We denote by $T'_F$ the sheaf of affine sections of $\hat T_F$.

In [T5], we gave a
  classification
  of affinely locally trivial
  affine bundles using Hochschild cohomology classes for
   a representation of
 $\pi_1(B)$,

 In this paper, we will give   another using Cech cohomological classes via gerbe theory.

In fact it seems to us that the second classification  fits best
to our problem. Inspired by the
philosophy of groupoid sheaves, we canonically associate to
any representation $\pi:\pi_1(B)\rightarrow Aff(F,\nabla_F)/T_F$,
 a gerbe with lien $T_F'$ which describes the
gluing problem related to the existence of affinely locally
trivial affine bundles associated to $\pi$.
   The $2-$cocycle described in
Giraud's  classification theorem of the associated gerbe,
 is given by an
element of $H^2(B,T_F')$. This class is the obstruction
to the existence of affinely locally trivial affine bundles associated
to $\pi$. When it vanishes, each element
of $H^1(B,T_F')$ defines an affinely locally trivial affine bundle.
In this case
 the classification of the affinely locally trivial affine bundles
is given by the orbits of  elements of $H^1(B,T_F')$ under a gauge
group.

After the classification of affine bundles, we classify
composition series of affine manifolds in which each map
$f_i:(M_{i+1},\nabla_{M_{i+1}})\rightarrow (M_i,\nabla_{M_i})$ is
an affinely locally trivial affine bundle. Since the
classification of affinely locally trivial affine bundle has been
done using gerbe theory, it is natural to think that the theory
involved in the classification of composition series of affine
manifolds is $n-$gerbe theory. In the last section of our work, we
build a commutative $n-$gerbe theory.

This is the plan of our paper:

$0.$ Introduction.

\medskip

I. AFFINE BUNDLES.

\medskip

$1.$ Background.

$2.$ Generality.

$3.$ The classification of affinely locally trivial affine bundles.

$4.$ The general case.

\medskip

II. COMPOSITION SERIES OF AFFINE MANIFOLDS.

\medskip

$1.$ $3$ composition series.

$2.$ The general case ($n-$composition series of affine
manifolds).

$3.$ The conceptualization (commutative $n-$gerbe theory).

\bigskip

{\bf I. AFFINE BUNDLES.}

\bigskip

{\bf 1.  Background.}

\bigskip

An $n-$ connected affine manifold $(M,\nabla_M)$, is an $n-$ connected
differentiable
manifold endowed with
a connection $\nabla_M$, whose curvature and torsion forms
vanish identically. The connection $\nabla_M$ defines on $M$ an atlas (affine)
 whose
transition functions are locally affine transformations of ${\R}^n$.

Let $(M,\nabla_M)$ and $(N,\nabla_N)$ be two affine manifolds respectively
associated to the affine atlas, $(U_i,\phi_i)$ and $(U'_j,\phi'_j)$.
An affine map between $(M,\nabla_M)$ and $(N,\nabla_N)$ is a differentiable
map $f:M\rightarrow N$ such that
$\phi'\circ
 f_{\mid{U_i}}\circ\phi_i^{-1}$
is an affine map. We denote by $App((M,\nabla_M),(N,\nabla_N))$ the set
of affine maps between $(M,\nabla_M)$ and $(N,\nabla_N)$, and by
$Aff(M,\nabla_M)$ the space of affine automorphisms of $(M,\nabla_M)$.

The affine structure of $M$ pulls back
 to its universal cover $\hat M$, and defines
on it an affine structure $(\hat M, \nabla_{\hat M})$, for which
the universal cover map $p_M:\hat M\rightarrow M$ is an affine map.
 The affine structure of
$(\hat M,\nabla_{\hat M})$ is defined by a local diffeomorphism
$D_M:\hat M\rightarrow {\R}^n$ called the developing map.

The developing map gives rise to a representation
$A_M:Aff(\hat M,\nabla_{\hat M})\rightarrow Aff({\R}^n)$ which makes
the following diagram commute

$$ \matrix
{(\hat M,\hat\nabla_{\hat M})&{\buildrel{g}\over{\longrightarrow}}&
(\hat M,\hat\nabla_{\hat M})\cr\downarrow D_M &\ \ \ \ & \downarrow
D_M\cr{\R}^n &{\buildrel{A_M(g)}\over{\longrightarrow}}&
{\R}^n} $$

where $g$ is an element of $Aff(\hat M,\nabla_{\hat M})$. The restriction of
$A_M$
to the fundamental group $\pi_1(M)$ of $M$, is the holonomy representation
$h_M$. The linear part $L(h_M)$ of $h_M$, is the linear holonomy of
$(M,\nabla_M)$. It is in fact the holonomy of the connection $\nabla_M$
in the classical sense.

\medskip

{\bf Definitions 1.1.}

- The affine manifold $(M,\nabla_M)$ is complete, if and only if the developing
map is a diffeomorphism. This is equivalent to saying that the connection
$\nabla_M$ is geodesically complete.

- The affine manifold $(M,\nabla_M)$ is unimodular, if its linear holonomy
lies in $Sl(n,{\R})$. Markus has conjectured
that a compact affine manifold is complete if and only if it is unimodular.

- Let $f$ and $g$ be two affine bundles with the same base space $(B,\nabla_B)$,
and respectively total spaces $(M,\nabla_M)$, and $(N,\nabla_N)$. An affine
isomorphism between $f$ and $g$, is an affine isomorphism
 between $(M,\nabla_M)$ and
$(N,\nabla_N)$ which sends a fiber of $f$ onto a fiber of $g$, and gives rise
to an automorphism of $(B,\nabla_B)$.

\bigskip

{\bf 2. Generalities.}

\medskip

This paragraph is devoted to some basic properties of affine bundles.
In the sequel, we will suppose that all the fibers of a given affine bundle
are  diffeomorphic to each other.

Let $f:(M,\nabla_M)\rightarrow (M',\nabla_{M'})$ be an affine map,
the map $f$ pulls back to a map $\hat f:(\hat M,\nabla_{\hat M})\rightarrow
(\hat M',\nabla_{\hat M'})$ which makes the following diagram
commute:

$$ \matrix
{(\hat M,\nabla_{\hat M})&{\buildrel{\hat f}\over{\longrightarrow}}&
(\hat M',\nabla_{\hat M'})\cr\downarrow p_M &\ \ \ \ & \downarrow
p_{M'}\cr (M,\nabla_M) &{\buildrel{f}\over{\longrightarrow}}&
(M',\nabla_{M'})}. $$

\medskip

{\bf Proposition 2.1. [T4].}
{\it Let $(M,\nabla_M)$ be the domain of an affine map.  Suppose that
$M$ is compact. We denote by $df_x$, the differential  $df$
of $f$, at $x$. Then the distribution $Df$ of $M$ defined by
$$
Df_x=\{v\in T_xM/ df_x(v)=0\}
$$
defines on $M$ an affine bundle whose fibers are the leaves of
the foliation defined  by $Df$.}

\medskip

{\bf Sketch of proof.}

As $M$ is compact, the space of fibers is a differentiable manifold, say $B$.
The transverse affine structure of the foliation
$Df$, pushes forward  to $B$ and defines on it
an affine connection $\nabla_B$, which makes the projection
$(M,\nabla_M)\rightarrow (B,\nabla_B)$ an affine map.

\medskip

 This  proposition implies that an affine bundle with
 compact total space gives rise to  a locally trivial differentiable
  bundle  by a
 well-known Ehresmann result [God].  Denote by $F$ the typical fiber. Applying the
  Serre bundle sequence to
 this bundle, we obtain the following short exact sequence:
 $$
 \pi_2(B)\rightarrow \pi_1(F)\rightarrow \pi_1(M)\rightarrow\pi_1(B)
 \rightarrow 1.
 $$
 If we suppose that $\pi_2(B)=1$, then we obtain that $\pi_1(M)$ is
 an extension of $\pi_1(B)$ by $\pi_1(F)$. In particular this
 happens when $(M,\nabla_M)$ is complete. Remark  that if $(M,\nabla_M)$
 is complete, then the fibers and the base of the induced bundle
  are also complete.

 Auslander has conjectured that the fundamental group of a compact and
 complete affine manifold is polycyclic. In [T4], we have conjectured
 that we can change the complete affine structure of a galoisian
 cyclic finite cover of a
  $n-$compact, $(n>2)$ and complete affine manifold to another complete one,
    so that it becomes  the domain of a non trivial affine map.
     Non trivial means that
 the distribution $Df$ is neither $0$, nor the whole space. This
 conjecture implies the Auslander conjecture (see [T4]). As mentioned
 in the introduction, the classification of affine compact bundles
 will conjecturally allow us to know  all compact and complete affine manifolds up to a finite cover.

 In fact, there are examples of affine manifolds, which are total spaces of
 more than one non isomorphic affine bundle. The following is an example
 of this situation in dimension $3$.

 Let $C=(e_1,e_2,e_3)$ be a basis of ${\R}^3$. Consider the subgroup
 $\Gamma$ of $Aff({\R}^3)$, generated by $f_1$, $f_2$ and $f_3$, whose
 expressions in $C$ are:

 $f_1(x,y,z)=(x+1,y,z)$

 $f_2(x,y,z)=(x,y+1,z)$

 $f_3(x,y,z)=(x+y,y,z+1)$.

 The quotient of ${\R}^3$ by $\Gamma$ is a compact affine manifold,
 $M^3$. The projections $p_2(x,y,z)=y$, and $p_3(x,y,z)=z$, define
   projections of $M^3$ over the circle endowed with its canonical complete
 structure. The bundles defined by those projections are not isomorphic.
 If they were isomorphic, there would exist an element of $Aff({\R}^3)$
 of the form $(x,y,z)\rightarrow (ax+b,f(y,z)+d(x))$ where $f$ is an element of
 $Aff({\R}^2)$, and $d$ a linear map ${\R}\rightarrow {\R}^2$
  which
  conjugates the map $(x,y,z)\rightarrow (x+1,y,z)$ to the map
 $(x,y,z)\rightarrow (x+1,y+z,z)$. This is evidently impossible.
 (For each bundle, we have adapted the expression of $\Gamma$ in a
 basis $(e_1',e_2',e_3')$ such that the vector subspace ${\R}e_1'$
  pulls forward on the
  base of each fibration).

 \medskip

 {\bf Definition 2.2.}

 Let $f:(M,\nabla_M)\rightarrow (B,\nabla_B)$ be an affine bundle,
 We will say that the bundle $f$ is an affinely locally trivial affine bundle, if and only
 if there exists an affine manifold $(F,\nabla_F)$ such that
  each element $x$ of $B$, is contained in an open set $U_x$, such that
  there exists an affine isomorphism
  $$f^{-1}(U_x)\rightarrow U_x\times (F,\nabla_F)$$
   and the
  restriction of the projection on $f^{-1}(U_x)$ is the first projection
  $U_x\times (F,\nabla_F)\rightarrow U_x$ via this isomorphism.

`\medskip

When the total space is compact, the last
definition is equivalent to saying  that one can build the Cech cocycle which defines the locally trivial
differentiable structure  of the affine bundle by affine maps.

In the previous examples the bundle defined by $p_3$ is  affinely locally
trivial, but not  the one defined by $p_2$.

Let $f:(M,\nabla_M)\rightarrow (M',\nabla_{M'})$ be an affine map,
where $M$ and $M'$ are respectively an $n$ and an $n'-$manifold.
The map $f$  pulls back to a map $\hat f:(\hat M,\nabla_{\hat M})\rightarrow
(\hat M',\nabla_{\hat M'})$. There exists an affine map
$f':{\R}^n\rightarrow {\R}^{n'}$
 which makes the following diagram
commute:

$$ \matrix
{(\hat M,\nabla_{\hat M})&{\buildrel{\hat f}\over{\longrightarrow}}&
(\hat M',\nabla_{\hat M'})\cr\downarrow D_M &\ \ \ \ & \downarrow
D_{M'}\cr{\R}^n &{\buildrel{f'}\over{\longrightarrow}}&
{\R}^{n'}}. $$

Let $(M,\nabla_M)$ be the total space of an affine bundle $f$.
The foliation $ {\cal F}_f$ defined by the leaves of $f$, pulls back to
 a foliation $\hat {\cal F}_f$ on $\hat M$,
which is the pull-back of a foliation $D_M(\hat {\cal F}_f)$ of ${\R}^n$,
 by parallel $l-$affine subspaces,
here $n$ and $l$ are the dimensions
of  $M$ and of the fibers of the bundle.

Write ${\R}^n={\R}^m\oplus{\R}^l$, where $m$ is the dimension
of the base of the bundle.

For every element $\gamma$ of $\pi_1(M)$, we  have:
$$
h_M(\gamma)(x,y)=(A_\gamma(x)+a_\gamma,B_\gamma(y)
+C_\gamma(x)+d_\gamma)
$$
as $D_M(\hat{\cal F}_f)$ is stable under the holonomy.

An element $\gamma$ of $\pi_1(M)$ which preserves a fiber of $\hat f$,
 preserves all the other fibers as the foliation ${\cal F}_f$ does not
  have holonomy.
 We obtain that
 $$
 h_M(\gamma)(x,y)=(x,B_\gamma(y)+C_\gamma(x)+d_\gamma).
 $$
 In fact we obtain a representation $\pi_1(F)\rightarrow Aff({\R}^l)$
 for each $x\in {R}^m$.

  We deduce that, as they have been supposed to be diffeomorphic,
    all the fibers have the same linear holonomy.

  The map
   $$
 \pi_1(F)\longrightarrow {\R}^l
 $$
$$
 \gamma\longrightarrow C_{\gamma}(x)+d_\gamma
 $$
 is a $1-$cocycle with respect to the linear holonomy.
 It defines the map

 $$
 r:{\R}^m\rightarrow H^1(\pi_1(F),{\R}^l)$$
 $$ x\rightarrow [C_\gamma(x)+d_\gamma]$$
 where $H^*(\pi_1(F),{\R}^l)$ is the $*$ cohomology  group,
 with  respect to the linear holonomy of the fibers.
 The cohomology class $r(x)$ is often called the radiance obstruction
 of the affine holonomy of the fiber over $p_B(x)$.

If the fibers are compact and complete,
 the image of $r$ is contained in the  algebraic subvariety  $L$, of
 $H^1(\pi_1(F),{\R}^l)$ defined by
 $$
 L=\{c\in H^1(\pi_1(F),{\R}^l)/ \Lambda^lc\neq 0\}.
 $$
 See [F-G-H] theorem 2.2.

 The fact that $f$ is an affinely locally trivial affine bundle, is equivalent
 to the fact that the map $r$ is a constant map when the total space is
 complete.

\medskip

{\bf Question.}
Are the fibers of an affine bundle isomorphic if its total space
is compact ?

\medskip

The following theorem was inspired by the last question:

\medskip

{\bf Theorem 2.3. [T5].}
{\it  Suppose that the total space of an affine bundle is
an $n-$compact and complete affine manifold, and
moreover
the fundamental group of the fibers are nilpotent.
 Then all the fibers are isomorphic  to each
other.}

\medskip

Given an element $\gamma$ of $\pi_1(B)$, and an element
$x\in {\R}^m$, the restriction  of the affine
holonomy representation of
$\pi_1(F)$ to $x\times {\R}^l$ and $(h_B(\gamma))(x))\times {\R}^l$
  are conjugated by an element
of $Aff({\R}^l)$, since they define the same affine
structure (we choose $x\in D_B(\hat B)$).
 This leads to define a gauge group for
the linear representation $L(h_F)$.

Consider the subgroup $G$ of automorphisms of $\pi_1(F)$ such that
for every element  $g$ of $G$, there is a linear map $B_g$ such that

$$
L(h_F)(g(\gamma))=B_g\circ L(h_F)(\gamma)\circ B_g^{-1}.
$$
The group which elements are $B_g$ will be called the gauge group of $L(h_F)$.

We associate to every $B_g\in G$ the following  linear map of
$H^p(\pi_1(F),{\R}^l)$: for each $c\in H^p(\pi_1(F),{\R}^l)$ we
define $B_g^*c$ by
$$
B_g^*c(\gamma_1,...,\gamma_p)=B_g(c(g^{-1}(\gamma_1),...,g^{-1}(\gamma_p))).
$$
Two  complete affine structures $\nabla_1$ and $\nabla_2$ on $F$, with same linear holonomy
$L(h_F)$, and holonomy
$h_1$ and $h_2$,
are isomorphic if and only if there is an element $B_g$ of the gauge group
such that $B_g^*c_1=c_2$, where $c_1$ and $c_2$ are respectively the
radiance obstruction of $h_1$ and $h_2$.

\medskip

{\bf 3. Affinely locally trivial affine bundles.}

\medskip

Recall that, an affine bundle with complete total space
 is said to be an affinely locally trivial
affine
bundle, if and only if its pull-back of the bundle to the universal
cover of the base is a trivial affine bundle.

In the sequel, $(M,\nabla_M)$ will be the compact   total space of an affinely locally trivial affine bundle, with
base space $(B,\nabla_B)$ and typical fiber $(F,\nabla_F)$.

The following proposition emphasizes the importance of the category of
affinely locally trivial affine bundles.

\medskip

{\bf Proposition 3.1. [T5]}
{\it Let $f$ be an affine bundle whose total space is a complete affine
$n-$manifold
(not necessarily compact). If we suppose that the
fibers are $2-$tori and moreover  their linear holonomy
is the linear holonomy of a complete structure of the $2-$torus distinct
from the flat Riemmannian one, then $f$ is an affinely locally  trivial
affine bundle.}

\medskip

Let's go back to the classification problem.

Recall that $\pi_1(F)$ is a normal subgroup of $\pi_1(M)$.
Let $\gamma$ and $\gamma_1$ be respectively two elements  of $\pi_1(F)$
and $\pi_1(M)$. We can write

$$
h_M(\gamma)=(x,B_\gamma(y)+d_\gamma)$$

and

$$
h_M(\gamma_1)(x,y)=(A_{\gamma_1}(x)+a_{\gamma_1},
B_{\gamma_1}(y)+C_{\gamma_1}(x)+d_{\gamma_1}),
$$
where $A_{\gamma_1}$ is an automorphism of ${\R}^m$, $B_{\gamma}$
and $B_{\gamma_1}$ are automorphisms of ${\R}^l$,
$C_{\gamma_1}:{\R}^m\rightarrow {\R}^l$ is a linear map,
$a_{\gamma_1}$ is an element of ${\R}^m$ and $d_\gamma$, $d_{\gamma_1}$
are elements of ${\R}^l$.

One has

$$
h_M(\gamma_1)^{-1}(x,y)=
(A^{-1}_{\gamma_1}(x)-A_{\gamma_1}^{-1}(a_{\gamma_1}),
B_{\gamma_1}^{-1}(y)-B_{\gamma_1}^{-1}(d_{\gamma_1})
-B_{\gamma_1}^{-1}C_{\gamma_1}(A^{-1}_{\gamma_1}(x)-
A^{-1}_{\gamma_1}(a_{\gamma_1}))).
$$

Using the fact that $h_M(\pi_1(F))$ is a normal subgroup of
$h_M(\pi_1(M))$, one obtains
$$
\matrix{h_M(\gamma_1^{-1})\circ h_M(\gamma)\circ h_M(\gamma_1)=\cr
(x,B_{\gamma_1}^{-1}B_\gamma B_{\gamma_1}(y)+
B_{\gamma_1}^{-1}B_\gamma C_{\gamma_1}(x)+
B_{\gamma_1}^{-1}B_\gamma(d_{\gamma_1})+B_{\gamma_1}^{-1}(d_\gamma)
-B_{\gamma_1}^{-1}C_{\gamma_1}(x)-B_{\gamma_1}^{-1}(d_{\gamma_1}))}.
$$

This implies that

$$
B_{\gamma_1}^{-1}B_\gamma C_{\gamma_1}(x)-
B_{\gamma_1}^{-1}C_{\gamma_1}(x)=0;
$$
 we deduce that $C_{\gamma_1}(x)\in H^0(\pi_1(F),{\R}^l)$.
 The linear space $Applin({\R}^m,H^0(\pi_1(F),{\R}^l))$ of linear
 maps ${\R}^m\rightarrow H^0(\pi_1(F),{\R}^l)$ has
 a left $\pi_1(B)-$module structure defined by
 $$
 \gamma_1'(D)=B_{\gamma_1}\circ D
 $$
 and a right $\pi_1(B)-$module structure defined by
 $$
 \gamma_1'(D)=D\circ A_{\gamma_1},
 $$
 where $\gamma_1$ is an element of $\pi_1(M)$ over an element
 $\gamma_1'$ of $\pi_1(B)$.

 \medskip

 We denote by $T_F$ the connected component of the group of affine maps
 of $(F,\nabla_F)$, which pull-back on translations of ${\R}^l$.
 The linear map of $H^0(\pi_1(F),{\R}^l)$ defined by
 $$
 t\rightarrow B_{\gamma_1}t
 $$
 induces a linear map of $T_F$. This induces a $\pi_1(B)$ left
 structure on $Applin({\R}^m,T_F)$. The right structure of
 $\pi_1(B)$ on $Applin({\R}^m,H^0(\pi_1(F),{\R}^l))$ also induces
 a right $\pi_1(B)$ structure on $Applin({\R}^m,T_F)$.
 Let ${\Z}\pi_1(B)$ be the group algebra of the group ${\pi_1(B)}$.
 The vector space $Applin({\R}^m,T_F)$ is endowed with a ${\Z}\pi_1(B)$
 Hochschild module structure.

The bundle is supposed to be affinely locally trivial; this implies
that his lifts on $\hat B$, is $\hat B\times (F,\nabla_F)$.
The action of $\pi_1(B)$ on $\hat B\times (F,\nabla_F)$ is made by
affine maps. We deduce that $\pi_1(F)$ is normal in $\pi_1(M)$,
and a representation $\pi_f:\pi_1(B)\rightarrow Aff(F,\nabla_F)/T_F$
induced by its action on $\hat B\times (F,\nabla_F)$.

\medskip

Recall the following problem stated in [Bry].
Let
$$
1\rightarrow A\rightarrow G\rightarrow H\rightarrow 1
$$
be an exact sequence of groups, where $A$
is commutative. Given an $H-$bundle over
the manifold $X$, we want to classify  all the  bundles over $X$, with
structural group $G$, which are lifts of the previous bundle.

\medskip

Our problem is quite similar: We have seen that
an affinely locally trivial affine bundle $f$ gives rise
to a representation
$\pi_f:\pi_1(B)\rightarrow Aff(F,\nabla_F)/T_F$.
We denote by  $\pi_f'$ the flat bundle  induced by $\pi_f$.
 Given such a representation $\pi:\pi_1(B)\rightarrow Aff(F,\nabla_F)/T_F$ with associated $Aff(F,\nabla_F)/T_F$
 bundle $\pi'$,
 we want to classify all affinely locally trivial affine
  bundles  associated.

   To make our theory fit in gerbe theory with
the  band $T_F'$, the sheaf of affine section of $\hat T_F$,
 we will consider first, the classification of
affine bundles up to $T_F'$ isomorphisms.

\medskip

Now let's  recall some general facts of sheaf of groupoids and descent
theory.  Our exposition follows the treatment of [Bry].
The general philosophy is to express  gluing conditions in terms of
covering maps. We remark that if a manifold is an affine manifold, so are its
covers,  the category of affine manifolds is  stable under affine fiber product.

For every locally isomorphic affine map
 $g:(Y,\nabla_Y)\rightarrow (B,\nabla_B)$, we can pull-back the
 affinely locally trivial affine bundle
 $\pi'$ over $(B,\nabla_B)$ in a bundle $\pi_{Y}'$ over $(Y,\nabla_Y)$
 (associated to $\pi_Y$).
  The total space of the pull-back is
 $Q= \pi'\times_B (Y,\nabla_Y)$ and the fiber map the canonical
 projection.

 Consider $Y\times_B Y$ with the two canonical projections
 $p_1,p_2:Y\times_B Y\rightarrow Y$.
 We denote by $p_i^*Q$ $i=\{1,2\}$ the pull-back of $Q$ by $p_i$.
 There is a natural isomorphism $\phi: p_1^*Q\rightarrow p_2^*Q$.
 This isomorphism satisfies the following cocycle condition

 $$
 p^*_{13}(\phi)=p^*_{23}(\phi)\circ p^*_{12}(\phi),\leqno (2)
 $$
 an equality of morphisms $p^*_1Q\rightarrow p_3^*Q$ of affine bundles
 over $Y\times_B Y\times_B Y$, where $p_1$, $p_2$ and $p_3$ are the canonical
 projections on the three factors, and $p_{12}$, $p_{13}$ and $p_{23}$
 are the canonical projections of $Y\times_B Y\times_B Y$ over $Y\times_B Y$.

Conversely, given an affine bundle $Q\rightarrow Y$ which satisfies
the condition $(2)$, we recover an affine bundle over $(B,\nabla_B)$.

In fact one obtains:

\medskip

{\bf Proposition 3.2.}
{\it Let $g:(Y,\nabla_Y)\rightarrow (B,\nabla_B)$ be a local isomorphism
 of affine manifolds. The pull-back functor $g^*$ induces an equivalence
of categories between the category of  affine bundles over $(B,\nabla_B)$
and the category of  affine bundle bundles over $(Y,\nabla_Y)$ equipped with a
descent isomorphism $\phi:p_1^*Q\rightarrow p_2^*Q$ satisfying the cocycle
condition (2).}

\medskip

 The general
definition of torsor adapted to this case is:

\medskip

{\bf D\'efinition 3.3.}
{\it A $T_F'$ torsor, will be a sheaf $H$ on $(B,\nabla_B)$,
 together with
a $T_F'$ action such that every point of $B$ has a neighborhood $U$
with the property that for every $V\subset U$ open,
the space $H(V)$ is an affine principal bundle with structural
group $T_F'{_{\mid  V}}$.}

\medskip

 The isomorphism classes of $T'_F$ torsors are given by $H^1(B,T_F')$,
  where
 $H^*(B,T_F')$ is the $*$ cohomology  group of $T_F'$ related to the usual Cech
 cohomology.

 As the sheaf $T_F'$ is a locally constant sheaf, the notion of $T_F'$ torsor
 in our case is similar to a notion of affine $T_F'$ bundle over
 $(B,\nabla_B)$.

 \medskip

We can associate to
a representation $\pi:\pi_1(B)\rightarrow Aff(F,\nabla_F)/T_F$,
 the following sheaf  of groupoids, $B_\pi$.
To every local affine isomorphism $(Y,\nabla_Y)\rightarrow (B,\nabla_B)$,
we associate the category $Y_\pi$, whose  objects are
affinely locally trivial affine bundles
 over $(Y,\nabla_Y)$ with typical fiber $(F,\nabla_F)$,
associated to $\pi_Y$. The
  (auto)morphisms are $T'_F-$automorphisms.

It is easy to show that the following properties are satisfied:

(i) For every diagram $(Z,\nabla_Z){\buildrel{g}\over{\longrightarrow}}
(Y,\nabla_Y){\buildrel{h}\over{\longrightarrow}}(B,\nabla_B)$ of local
affine isomorphisms, there is a functor
$g^{-1}:Y_\pi\rightarrow Z_\pi$;

(ii)
For every diagram
 $(W,\nabla_W){\buildrel{k}\over{\longrightarrow}}
(Z,\nabla_Z){\buildrel{g}\over{\longrightarrow}}
(Y,\nabla_Y){\buildrel{h}\over{\longrightarrow}}(B,\nabla_B)$
there is an invertible natural transformation
 $\theta_{g,k}: k^{-1}g^{-1}\rightarrow (gk)^{-1}$.

 \medskip

 This makes our category a presheaf category. Moreover
  properties a la $(2)$ are satisfied to ensure that some kind
   of Haefliger $1-$cocycles
  is satisfied, in order to make our presheaf of category a sheaf of
  category.

  One can more usually define a sheaf of category. It is a map
   $C$ on the family of open subsets of $B$
  $$
  U\longrightarrow C(U)
  $$
   which assigns to any open subset $U$ of $B$ a category $C(U)$.

   For every open subset $V\subset U$, there is a composition of
   morphisms from $C(U)$ to $C(V)$.
   When $U=V$, this composition is just the compostion of morphisms.
   This defines the presheaf of category.
   Moreover a descent condition is needed to make the presheaf a sheaf.

  \medskip

 In fact, our sheaf of category is a gerbe with band $T'_F$. It means that
 the following properties are satisfied

 (G1) Given any object of $Y_\pi$, the sheaf of local automorphisms of this
 object is a sheaf of groups which is locally isomorphic to $T_F'$.

 (G2) Given two objects $Q_1$ and $Q_2$ of $Y_\pi$, there exists a
 local isomorphism surjective map $g:Z\rightarrow Y$ such that $g^{-1}Q_1$ and
 $g^{-1}Q_2$ are locally isomorphic.

 (G3) There is a  local isomorphism surjective affine map $Y\rightarrow X$
 such that the category $Y_\pi$ is not empty.

One say that our sheaf of category is a gerbe with band or lien $T_F'$.

 \medskip

{\bf Remark.}

To ensure the axiom (G3) to be satisfied, one may show an affinely
local trivial affine bundle with typical fiber $(F,\nabla_F)$,
 over $\hat B$. This bundle is just the trivial one.

 Let's now state the first classification theorem which
 is an adaptation of the Giraud
 classification theorem.

\medskip

 {\bf Theorem 3.4.}
 {\it
 The set of equivalence classes of the gerbes  is in one to one  correspondence
 with
 $H^2(B,T_F')$.}

\medskip

{\bf Proof.}

Let's consider a cover $(U_i)_{i\in I}$ of $B$ by open $1-$connected
affine charts. The $T'_F-$automorphisms of
an object  $P_i$ of
$C(U_i)$, is isomorphic to the restriction of $T'_F$ to $U_i$.

There is a $T'_F-$isomorphism
$$u_{ij}:(P_i)_{\mid C( U_{ij})}\rightarrow
( P_j)_{\mid C(U_{ij})}$$
in the category $C(U_{ij})$. We define a section
$h_{ijk}$ of $T'_F$ by,
$$
h_{ijk}=u^{-1}_{ik}u_{ij}u_{jk}\in Aut(P_k)\simeq T_F'.
$$

In fact $h=(h_{ijk})$ is a $T_F'-$valued Cech $2-$cocycle.
The
corresponding class in $H^2(B,T_F')$ is independent of  all the choices.
We will show that this correspondence defines
 an isomorphism between the group of equivalence classes,
and the set of isomorphic  gerbes with band $T_F'$.

To show the injectivity of this map, one remarks that if the cohomology
class defined by $h$ is trivial, then one can modify the isomorphisms
$u_{ij}$ such that $u_{ik}=u_{ij}u_{jk}$. We then obtain an affine $T_F'$ torsor
over $(B,\nabla_B)$ which represents a trivial gerbe.

To prove the surjectivity, we construct a gerbe associated to a $2-$Cech
cocycle $h=(h_{ijk})$ with values in $T_F'$. It   is sufficient to find
a family of elements $u_{ij}$ of $T_F'-$automorphisms of $U_{ij}$ such that
the condition $u_{ik}^{-1}u_{ij}u_{jk}=h_{ijk}$ is satisfied.

\medskip

This is our classification theorem for affinely locally trivial affine bundles,
up to $T_F'-$isomorphisms.

\medskip

{\bf Theorem 3.5.}
{\it Let
 $\pi:\pi_1(B)\rightarrow Aff(F,\nabla_F)/T_F$ be a representation. Then
 there are affine bundles over $\pi'$,
  if and only if its associated gerbe  is trivial. In this case the
  $T'_F-$isomorphism classes of affine bundles are given by
  $H^1(B,T_F')$.}

 \medskip

 {\bf Proof.}

If there exists an affine bundle associated to the representation $\pi$,
 the boundary of the cocycle which defines the fibration represents the associated
 gerbe, so this gerbe is trivial.

 On the other hand, the $2-$cocycle associated can be described as
 follows:

 Consider a trivialisation of the flat bundle  $\pi'$, associated to
 $\pi$.

 For every $i,j$ such that $U_i\cap U_j\neq\emptyset$, we have
  $$U_i\cap U_j\times Aff(F,\nabla_F)/T_F\longrightarrow
   U_i\cap U_j\times Aff(F,\nabla_F)/T_F$$
  $$(x,y)\longrightarrow (x,g'_{ij}(y)).$$
  We denote by $g_{ij}(x)$, an element of $Aff(F,\nabla_F)$ over
  $g'_{ij}$ which depends affinely on $x$.
 We set
 $$
 h_{ijk}=g_{ik}^{-1}g_{ij}g_{jk}.
 $$

 We have seen that if the cocycle is trivial, one can find a family
 of maps $w_{ij}:U_i\cap U_j\rightarrow T'_F$ such that
 $$
 (g_{ik}+w_{ik})=(g_{ij}+w_{ij})(g_{jk}+w_{jk}).
 $$
 Consider the
family of maps
 $$
\phi_{ij}: U_i\cap U_j\times (F,\nabla_F)\longrightarrow
U_i\cap U_j\times (F,\nabla_F)$$

 $$
 (x,y)\longrightarrow (x, (g_{ij}+w_{ij}(x))(y))
 $$

We have
$$
\phi_{ij}\circ\phi_{jk}(x,y)=(x,(g_{ij}g_{jk}+w_{ij}(x)+g_{i}w_{jk}(x))(y)).
$$

 One sees that the Cech cocycle condition is  verified.

For every other $1-$cocycle $(v_{ij})$, one gets an affine bundle
by setting
$$
\phi'_{ij}:U_i\cap U_j\times (F,\nabla_F)
\longrightarrow U_i\cap U_j\times (F,\nabla_F)
$$
$$
(x,y)\longrightarrow (x, (g_i+(w_{ij}+v_{ij})(x))(y)).
$$

\medskip

Two different cocycles used to define  affinely locally trivial
 affine bundles can define isomorphic affine bundles.

 Let $f_1$ and $f_2$ be two isomorphic affinely locally trivial
 affine bundles whose total
 spaces are $n-$compact, and  which induce the same sheaf $T_F'$.
 There is an affine transformation $g$ of $\hat f_1$ which preserves the
 foliation $\hat {\cal F}_{f_1}$ (where $\hat{\cal F}_{f_1}$
  is the pull-back of the
  the foliation induced by $f_1$) and conjugates
 the Deck transformations which defines the
 total space of $f_1$, in those which
 define the one of  $f_2$. As  both bundles induce
 $T_F'$, their induced representions
 $\pi_1(B)\rightarrow Aff(F,\nabla_F)/T_F$ coincide. The cohomology class
 defined by $f_1$ is changed in the class defined by $f_2$ by $g$, one has

 \medskip

 {\bf Proposition 3.6.}
{\it   The isomorphim
 classes of affine bundles are given by the quotient of $H^1(B,T_F')$
 by the action of a gauge group. This group is the group
  of affine automorphism of $\hat f_1$ which preserve its fibers,
  are pulls back of automorphisms of $(B,\nabla_B)$ and give rise to
 the same bundle $\pi'$.}

\medskip

{\bf 4. The general case.}

\medskip

In the previous section of our paper, we have classified affinely locally
trivial affine bundles. In the following section, we consider
the more general
 situation, when the total space is supposed to be only compact and
complete.

Given an affine bundle with compact and complete total space say $(M,\nabla_M)$ and
base space $(B,\nabla_B)$, we have seen that the fibers inherit
affine structures from the total space which are not necessarily isomorphic,
but which have the same linear holonomy. Let $\gamma$ be an element of
$\pi_1(M)$, set ${\R}^n={\R}^m\oplus{\R}^l$ where $n$, $m$, and $l$ are
respectively the dimension of $M$, $B$ and the typical fiber $F$.
We have:

$$
h_M(\gamma)(x,y)=(A_{\gamma}(x)+a_\gamma,B_\gamma(y)+C_\gamma(x)+d_\gamma),
$$
where $A_\gamma$ and $B_\gamma$ are respectively affine automorphisms
of ${\R}^m$ and ${\R}^l$, $a_\gamma$ and $d_\gamma$ are respectively
elements of ${\R}^m$ and ${\R}^l$, and $C:{\R}^m\rightarrow {\R}^l$ is
a linear map.

If $\gamma$ lies in $\pi_1(F)$, then $A_\gamma=I_{{\R}^m}$ and $a_\gamma=0$.

Now consider an
$m-$affine
manifold $(B,\nabla_B)$ compact and complete,
a compact $l-$manifold $F$ and a
 representation
$L(h_F):\pi_1(F)\rightarrow Gl(l,{\R})$ which is the
linear holonomy
of a complete affine structure of $F$.
We want to classify all affine bundles
with complete total space
 with base $(B,\nabla_B)$ and whose
fibers are diffeomorphic to $F$ and inherit affine structures from
the total space whose linear holonomy is $L(h_F)$.

 The natural question which arises is:
 Does every map $r:{\R}^m\rightarrow H^1(\pi_1(F),{\R}^l)$ give rise to an affine
 bundle ?

 For every element $\gamma_1$ of  $\pi_1(B)$,
 the affine representations defined by the cocycles $r(x)$ and $r(\gamma_1(x))$ must be
 isomorphic.

 Recall that we have defined a gauge group $G$ of the representation
 $L(h_F)$ as follows: it is a group of linear maps such
 that for every element $B_g\in G$ there is an automorphism $g$ of $\pi_1(F)$
   which satisfies
  $$
  L(h_F)(g(\gamma))=B_gL(h_F)(\gamma)B_g^{-1}.
  $$
The group $G$ acts on $L(h_F)$ one cocycles by setting
$$
(B_g^*c)(\gamma)=B_gc(g^{-1}(\gamma)).
$$
Let $\gamma_1$ and $\gamma$ be respectively  elements of $\pi_1(M)$ and
$\pi_1(F)$. We have
$$
\matrix{
\gamma_1\circ\gamma\circ\gamma_1^{-1}(x,y)=\cr
(x,B_{\gamma_1}B_\gamma B_{\gamma_1}^{-1}(y)-
B_{\gamma_1}B_\gamma B_{\gamma_1}^{-1}
C_{\gamma_1}(A_{\gamma_1}^{-1}(x)-A_{\gamma_1}^{-1}(a_{\gamma_1}))
-B_{\gamma_1}B_\gamma B_{\gamma_1}^{-1}(d_{\gamma_1})+\cr
B_{\gamma_1}C_\gamma(A_{\gamma_1}^{-1}(x)-A_{\gamma_1}^{-1}(a_{\gamma_1}))+
B_{\gamma_1}(d_\gamma)
+C_{\gamma_1}(A_{\gamma_1}^{-1}(x)-A_{\gamma_1}^{-1}(a_{\gamma_1}))+
d_{\gamma_1})}.
$$
The map
$$i'({\gamma_1}):\pi_1(F)\longrightarrow \pi_1(F)
$$

$$
\gamma\longrightarrow \gamma_1\gamma\gamma_1^{-1}
$$
is an automorphism    associated to the element of the gauge group
 $B_{\gamma_1}$. If $\gamma_1$ lies in $\pi_1(F)$ the induced map
on $H^1(\pi_1(F),{\R}^l)$ is trivial. We deduce a map
$$i:\pi_1(B)\longrightarrow Gl(H^1(\pi_1(F),{\R}^l))$$
$$\gamma_1\longrightarrow ai'(\gamma_1'),$$ where $\gamma_1$ pulls back
to $\gamma_1'$ and $ai'(\gamma_1')$ is the action of $\gamma_1'$ on
$H^1(\pi_1(F),{\R}^l)$ induced by $B_{\gamma_1'}$.

We have $r(A_{\gamma_1'}(x)+a_{\gamma_1'})=[B_{\gamma_1'}r(x)].$

It follows that the following square is commutative

$$ \matrix
{{\R}^m&{\buildrel{\gamma_1}\over{\longrightarrow}}&
{\R}^m\cr\downarrow r &\ \ \ \ & \downarrow
r\cr H^1(\pi_1(F),{\R}^l) &{\buildrel{i(\gamma_1)}\over{\longrightarrow}}&
H^1(\pi_1(F),{\R}^l)}.$$

The representation $i$ allows also to construct a  bundle over
$(B,\nabla_B)$ with typical fiber $H^1(\pi_1(F),{\R}^l)$.
We also denote by $i$ this bundle. The map
$r$ can also be viewed as a section of this bundle.

\medskip

We will classify all affine bundles for a
given  representation $i$, and a map $r$ such that each $r(x)$ defines
an affine structure.

The map $r$ defines a representation $\pi_1(F)\longrightarrow Aff({\R}^n)$
such that the quotient ${\R}^n/\pi_1(F)$ is an affine bundle over
${\R}^m$.
We assume that for each $\gamma\in \pi_1(B)$, there is an element
of $Aff({\R}^n/\pi_1(F))$
 which induces $i(\gamma)$.

Let $U$ be an open set of $(B,\nabla_B)$. We define
the following sheaf of categories
$$
U\longrightarrow C(U).
$$
Where $C(U)$ is the set of affine bundles such that the canonical bundle
with typical fiber $H^1(\pi_1(F),{\R}^l)$ associated, is the restriction
of $i$ to $U$, and whose lifts on the universal cover of $U$ is the restriction
of ${\R}^n/\pi_1(F)$ to it.

The sheaf of categories $C$ is a gerbe with band $Aff({\R}^n/\pi_1(F))_0$,
which denotes the sheaf induced by the
 connected component of the affine automorphisms of
${\R}^n/\pi_1(F)$ which pushes forward on the identity of ${\R}^m$.

Let us explain why the band is $Aff({\R}^n/\pi_1(F))_0$.

Consider a trivialization of $i$,
$$
\bar h_{ij}:U_i\cap U_j\times H^1(\pi_1(F),{\R}^l)\longrightarrow
U_i\cap U_j\times H^1(\pi_1(F),{\R}^l)
$$
$$
(x,y)\longrightarrow (x,\bar h_{ij}(y)).
$$
Here the $U_i$ are connected open sets of affine charts. So we can
restrict the bundle ${\R}^n/\pi_1(F)$ to each $U_i$. We denote by
$U_i^F$ this restriction.

To $\bar h_{ij}$ we associate an element $h_{ij}$ of $Aff(U_i^F)$ which
pushes forward on the identity of $U_i$, and gives rise to $\bar h_{ij}$.
The maps
$$ h_{ik}^{-1}h_{ij}h_{jk}$$ are the obstructions of the existence of
an affine bundle associated to $i$ and $r$.

\medskip

{\bf Proposition 4.1.}
{\it The map $h_{ijk}=h_{ik}^{-1}h_{ij}h_{jk}$ is an element of the restriction
of $Aff({\R}^n/\pi_1(F))_0$ to $U_i\cap U_j\cap U_k$.}

\medskip

{\bf Proof.}

We deduce this, from the fact that,
the map $h_{ijk}$ gives rise to the identity of $H^1(\pi_1(F),{\R}^l)$, since
it is shown in [T2] that the set of affine automorphisms which commutes with the
holonomy of a compact and complete affine manifold is a cover of the connected
component of its affine automorphism group.

\medskip

This is our classification theorem in the general case.

\medskip

{\bf Theorem 4.2.}
{\it For each representation $i$ and map $r$, there is a $2-$Cech
cocycle which is the obstruction of the existence of an affine bundle
associated to $i$ and $r$. When it vanishes the set of isomorphisms
classes of affine bundles are given by the orbit of elememt
of $H^1(B,Aff({\R}^n/\pi_1(F))_0)$ under a gauge group.}

\bigskip

{II. COMPOSITION SERIES OF AFFINE MANIFOLDS.}

 \bigskip

 Recall that a $n-$affine manifold is said to be complete if and only if
 the connection $\nabla_M$ is complete. This is equivalent
 to saying  that $M$ is the quotient
 of ${\R}^n$ by a group $\Gamma_M$
 of affine automorphisms which act properly and freely
 on ${\R}^n$. In this part, we will only consider complete affine manifolds.

 The representation $h_M:\Gamma_M=\pi_1(M)\rightarrow Aff({\R}^n)$ is called
 the holonomy of the affine manifold $(M,\nabla_M)$. Its
  linear part $L(h_M)$ is called
 the linear holonomy.

 It has been conjectured by Auslander that the fundamental group of a compact
 and complete affine manifold is polycyclic.

 In [T4], we have conjectured that each compact and complete $n>2$ affine
 manifold $(M,\nabla_M)$ has a finite cyclic and galoisian cover $M'$ endowed
 with a complete affine structure $(M',\nabla'_{M'})$
 eventually different from the pull back such that $(M',\nabla_{M'})$ is
 the source space of a non trivial affine map.
 Non trivial means that the fibers of $f$ are neither $M'$ nor points of $M'$.
 This conjecture implies the Auslander conjecture.

 If we suppose that the source space $(M,\nabla_M)$,
 of a non trivial affine map is compact, then $(M,\nabla_M)$ is the
 source space of a non trivial affine surjection $f$ over a manifold $(B,\nabla_B)$.
 We deduce from a well-known Ehresmann theorem, that $f$ is also a locally
 trivial differentiable fibration. All the fibers of an affine fibration
 inherit affine structures from $(M,\nabla_M)$ with same linear holonomy.

 The last conjecture leads to the following problem:
 Given two affine manifolds $(B,\nabla_B)$
 and $(F,\nabla_F)$ classify
every affine surjection $f:(M,\nabla_M)\rightarrow (B,\nabla_B)$
such that the differentiable structure of the fibers is $F$ and
their linear holonomy is the one of $(F,\nabla_F)$. Or more
generally, Given $n$ affine manifolds $(F_i,\nabla_{F_i})$,
classify all
 composition series $(M_{n+1},\nabla_{M_{n+1}})\rightarrow (M_n,\nabla_{M_n})
 \rightarrow...\rightarrow (M_1,\nabla_{M_1})$
 such that every map $f_i:(M_{i+1},\nabla_{M_{i+1}})\rightarrow (M_i,\nabla_{M_i})$
 in the last sequence is an affine surjection with fiber diffeomorphic to
 $F_i$, and which inherits from $(M_{i+1},\nabla_{M_{i+1}})$ an affine structure
 wich linear holonomy  is the linear holonomy of $(F_i,\nabla_{F_i})$.

 We have classify in the first part affine surjections with compact total spaces using gerbe
theory. The purpose of this part is to classify affine
 compostion series of affine manifolds.
We restrict to composition series $(M_n,\nabla_{M_n})\rightarrow
(M_{n-1},\nabla_{M_{n-1}})\rightarrow... \rightarrow
(M_1,\nabla_{M_1})$ such that the projection
$f_i:(M_{i+1},\nabla_{M_{i+1}})\rightarrow (M_i,\nabla_{M_i})$ is
an affinely locally trivial affine (a.l.t) fibration. This means
that the holonomy of the fibers is fixed.

Two composition series
$$({M_n}^j,\nabla_{{M_n}^j})\rightarrow
(M_{n-1}^j,\nabla_{M_{{n-1}^j}})\rightarrow...(
{M_2}^j,\nabla_{{M_2}^j})\rightarrow (M_1,\nabla_{M_1}), j=1,2$$
are equivalent if and only if the bundle ${f_i}^1$ and ${f_i}^2$
are isomorphic in respect  to $T_{F_{i}}$ isomorphisms. This means
that we consider isomorphisms of affine fibrations which act by
translations on the fibers and project on the identity of the base
space.

\medskip

  The
classification of affinely locally trivial affine bundles where
made using commutative gerbe theory. It is natural to think that
the classification of affinely locally $n-$ composition series
must be done using $n-$gerbe theory. In this part, we will give a
classification of $n-$ composition series of affine manifolds, and
after conceptualize the ideas involved to give a theory of
commutative $n-$gerbes.

\bigskip

\centerline{\bf 1. $3$ compostion series of affine manifolds.}

\bigskip

Let's recall first the classification  of affinely locally trivial
affine bundles up to translational isomorphisms with given fiber
and base space.

We have two affine manifolds $(B,\nabla_B)$ and $(F,\nabla_F)$
which represent respectively the base space and the fiber of the
affine bundles we intend to classify.

The structure of a locally trivial affine bundle gives rise to a
map $\pi_{BF}:\pi_1(B)\rightarrow Aff(F,\nabla_F)/T_F$ which
defines a locally flat bundle $PB_F$ over $B$ with typical fiber
$Aff(F,\nabla_F)/T_F$. This principal bundle induces a flat bundle
$BT_F$ over $B$ with typical $T_F$.

Let $T'_F$ be the sheaf of affine sections of $BT_F$. We define a
commutative gerbe  $C$ with lien $T'_F$ as follows:

To each open set $U$ of $B$, we associate the category $C(U)$ of
affinely locally trivial affine bundles with typical fiber
$(F,\nabla_F)$ such that the canonical $UT_F$ bundle associated to
it, is the restriction of $BT_F$ to $U$.

To represent the classifying two cocycle associated to $C$, we
consider an open covering $U_k$ of $B$ by  connected affine
charts, in each category $C(U_k)$ we choose an objet which is an
affine bundle isomorphic to $U_k\times (F,\nabla_F)$.

The trivialization of the bundle $PB_F$ gives rise to:
$$
U_k\cap U_l\times Aff(F,\nabla_F)/T_F\longrightarrow U_k\cap
U_l\times Aff(F,\nabla_F)/T_F
$$
$$
(x,y)\longmapsto (x,\bar t_{kl}(y)),
$$
where $\bar t_{kl}$ is an element of $Aff(F,\nabla_F)/T_F$. Taking
for each $k,l$ a map $t_{kl}$ over $\bar t_{kl}$ which depends
affinely of $x$, we obtain:

$$
U_k\cap U_l\times (F,\nabla_F)\longrightarrow U_k\cap U_l\times
(F,\nabla_F)
$$
$$
(x,y)\longmapsto (x,t_{kl}(x)y)
$$
Then we have the family of maps
$$
t_{klm}:U_k\cap U_l\cap U_m\longrightarrow T_F
$$
$$
x\longmapsto t_{kl}t_{lm}t_{km}^{-1}
$$
This family of maps $t_{klm}$ is a $2-$cocycle which classifies
the gerbe $C$. It is the obstruction to the existence of a locally
trivial affine bundle over $(B,\nabla_B)$ associated to $BT_F$. In
this case, the set of isomorphic classes of translational affine
bundles (or the classes of $T_F$ isomorphic bundles) with typical
fiber $(F,\nabla_F)$ and base space $(B,\nabla_B)$ are given by
the Cech cohomology group $H^1(B,T'_F)$ of the sheaf $T'_F$.

\bigskip

{\bf Remark.}

Consider the composition serie $(M_3,\nabla_{M_3})\rightarrow
(M_2,\nabla_{M_2})\rightarrow (M_1,\nabla_{M_1})$. The fact that
the affine bundles $(M_3,\nabla_{M_3})\rightarrow
(M_2,\nabla_{M_2}) $ and the  affine bundle
$(M_2,\nabla_{M_2})\rightarrow (M_1,\nabla_{M_1})$ are affinely
locally trivial affine bundles does not imply that the bundle
$(M_3,\nabla_{M_3})\rightarrow (M_1,\nabla_{M_1})$ is a locally
trivial affine bundle.

This can be illustrated by the following example. Consider the
subgroup $\Gamma$ of ${\R}^3$ generated by the three maps $f_1$,
$f_2$ and $f_3$ defined by
$$
f_1(x,y,z)=(x+1,y,z)
$$
$$
f_2(x,y,z)=(x,y+1,z)
$$
$$
f_3(x,y,z)=(x+y,y,z+1).
$$ In the canonical basis $(e_1,e_2,e_3)$ of ${\R}^3$.
The quotient of ${\R}^3$ by $\Gamma$ is a three compact affine
manifold $M^3$. The projection $p_1$ of ${\R}^3$ on its subvector
space  $V_2$ generated by $e_2$ and $e_3$ parallel to the one
generated by $e_1$, defines an affinely locally trivial affine
bundle over the torus $T_2$ endowed with its canonical riemannian
flat structure.

The projection $p_2$ of $V_2$ on the line generated by $e_2$
parallel to the one generated by $e_3$ defines an affinely locally
trivial affine bundle of the torus over the circle.

It is easy to see that the projection $p_2\circ p_1$ defines an
affine bundle which is not an affinely locally trivial affine
bundle over the circle.

\bigskip

Before to go to the general case, we will treat $3-$ series of
composition. So we have a sequence $(M_3,\nabla_{M_3})\rightarrow
(M_2,\nabla_{M_2})\rightarrow (M_1,\nabla_{M_1})$ of affine maps
such that $f_{i-1}:(M_i,\nabla_{M_i})\rightarrow
(M_{i-1},\nabla_{M_{i-1}})$ defines an affinely locally trivial
affine bundle.

We have supposed that the total spaces of our bundles are compact
and complete affine manifolds. This implies that $\pi_1(F_2)$ is
normal in $\pi_1(M_3)$ and $\pi_1(F_1)$ is normal in $\pi_1(M_2)$,
thus we have the following exact sequences

$$
1\rightarrow \pi_1(F_1)\rightarrow\pi_1(M_2)\rightarrow
\pi_1(M_1)\rightarrow 1,
$$
and
$$
1\rightarrow \pi_1(F_2)\rightarrow
\pi_1(M_3)\rightarrow\pi_1(M_2)\rightarrow 1.
$$

 We denote by $n_i, i=1,2,3$ the dimensions of $M_i$, and by
$l_i, i=1,2$ the dimensions of $F_i$. We put
${\R}^{n_3}={\R}^{n_1}\oplus{\R}^{l_1}\oplus{\R}^{l_2}$,

Let $\gamma$ be an element of  $\pi_1(M_3)$;  the
$(M_3,\nabla_{M_3})$ holonomy's
 action of $\gamma$ on ${\R}^{n_3}$
is given by
$$
\gamma(x,y,z)=(A_1^{\gamma}(x_1)+a_1^{\gamma},A_2^{\gamma}(x_2)+B_2^{\gamma}
(x_1)+a_2^{\gamma},A_3^{\gamma}(x_3)+B_3^{\gamma}(x_1,x_2)+a_3^{\gamma}),
$$
where $A_1^{\gamma}$ is an automorphism of ${\R}^{n_1}$,
$A_i^{\gamma}$, $i=2,3$ is an automorphism of ${\R}^{l_i}$,
$i=2,3$, $B_2^{\gamma}:{\R}^{n_1}\rightarrow {\R}^{l_1}$ is a
linear map, and $B_3^{\gamma}:{\R}^{n_1}\oplus
{\R}^{l_1}\rightarrow {\R}^{l_2}$ is a linear map.

If $\gamma$ belongs to $\pi_1(F_2)$, then the restriction of
$\gamma$ to ${\R}^{n_1}\oplus{\R}^{l_1}$ is the identity and
$B_3^{\gamma}=0$, since we have supposed the fibration
$(M_3,\nabla_{M_3})\rightarrow (M_2,\nabla_{M_2})$ to be an
affinely locally trivial affine fibration.

 \medskip

 The  holonomy of $(M_2,\nabla_{M_2})$ is given by the action of $\pi_1(M_2)$
 on ${\R}^{n_1}\oplus{\R}^{l_1}$ induced by the holonomy representation of
 $(M_3,\nabla_{M_3})$. In fact the restriction of $h_{M_3}$ to
 ${\R}^{n_1}\oplus {\R}^{l_1}$ factor through $\pi_1(M_2)$.

 We deduce that the $(M_2,\nabla_{M_2})$ holonomy's action of an
  element $\gamma$ of $\pi_1(M_2)$,  is

$$
(A_1^{\gamma}(x)+a_1^{\gamma},A_2^{\gamma}(y)+B_2^{\gamma}(x)+d_2^{\gamma}).
$$
If $\gamma$ is an element of $\pi_1(F_1)$, then $A_1^{\gamma}=id$,
$a_1^{\gamma}=0$ and $B_2^{\gamma}=0$ since we have supposed that
the fibration $(M_2,\nabla_{M_2})\rightarrow (M_1,\nabla_{M_1})$
is an affinely locally trivial affine fibration.

\medskip

Now, considering the holonomy of $(M_2,\nabla_{M_2})$, we write
the fact that $\pi_1(F_1)$ is normal in $\pi_1(M_2)$, then we
obtain that the image of $B_2^{\gamma}$ is contained in
$H^0(\pi_1(F_1),{\R}^{l_1})$. Here the group
$H^*(\pi_1(F),{\R}^{l_1})$, is the cohomology group of
$\pi_1(F_1)$ related to its linear holonomy.

Considering the holonomy of $(M_3,\nabla_{M_3})$, and writing that
$\pi_1(M_2)$ is a normal subgroup of $\pi_1(M_3)$,  we obtain that
$B_3^{\gamma}$ is contained in $H^0(\pi_1(F_2),{\R}^{l_2})$.

We have the representation
$$
\pi_1(M_1)\longrightarrow Aff(F_1,\nabla_{F_1})/T_{F_1}
$$ given by the first fibration. It leads to a flat bundle
$\pi_{11}(M_1)$ over $M_1$ with typical fiber
$Aff(F_1,\nabla_{F_1})/T_{F_1}$.

The bundle $(M_3,\nabla_{M_3})\rightarrow (M_2,\nabla_{M_2})$,
gives rise to a representation
$$\pi_1(M_2)\rightarrow Aff(F_2,\nabla_{F_2})/T_{F_2}$$ this leads to a
representation
$$\pi_{F_1}:\pi_1(F_1)\rightarrow Aff(F_2,\nabla_{F_2})/T_{F_2}$$
and to a flat principal bundle $\pi_{12}(F_1)$ over $F_1$ with
typical fiber $Aff(F_2,\nabla_{F_2})/T_{F_2}$ associated to
$\pi_{F_1}$. Since $\pi_1(F_1)$ is a normal subgroup of
$\pi_1(M_2)$, the action of $\pi_1(M_1)$ on this last bundle leads
to a representation
$$\pi_1(M_1)\rightarrow Aut(\pi_{12}(F_1))$$
where $Aut(\pi_{12}(F_1))$ is the group of automorphisms of the
bundle $\pi_{12}(F_1)$.

\bigskip

We will classify $3-$composition series of affine manifolds given:

The base space $(M_1,\nabla_{M_1})$ and the fibers spaces
$(F_1,\nabla_{F_1})$, $(F_2,\nabla_{F_2})$.

The representations $\pi_1(M_1)\rightarrow
Aff(F_1,\nabla_{F_1})/T_{F_1}$, $\pi_1(F_1)\rightarrow
Aff(F_2,\nabla_{F_2})/T_{F_2}$, and the representation
$\pi_1(M_1)\rightarrow Aut(\pi_{12}(F_1))$.

\bigskip

The main tool we will use to make the classification of $3$
composition series is $2$ gerbe, on this purpose, let recall some
definitions

\medskip

{\bf Definitions 1.1.}

A   category is  a  family of objects  and for each pair of
objects $X,Y$, a set of arrows $Hom(X,Y)$, which satisfy usual
rules.

An $2$  category is of a family of objects, and for each pair of
objects $X,Y$ the arrows is a  category $C(X,Y)$  which satisfies
usual rules see [Bry-Mc].

An $2$ gerbe on a manifold $M$ is a sheaf of $2$ categories  $C$
on $M$ which satisfies the following see [Bry-Mc], [Bre]:

($2.1$)For  every element $z$ of $M$, there exists
 an open set $U_z$ which contains
$z$ such that $C(U_z)$ is not empty.

($2.2$)Let $U$ be an open set, for any pair of object $x$ and $y$,
contained in $C(U)$, there is an open covering $(U_i)_{i\in I}$ of
$U$ such that, for any $i$, the set of arrows between the
restriction of $x$ and $y$ to $U_i$ is not empty.

($2.3$) For any $1$ arrow $f: x\rightarrow y$ in $C(U)$, there is
an inverse $g:y\rightarrow x$ up to a $2$ arrow.

($2.4$)The two arrows are invertible.

\medskip

We will associate to our problem a sheaf of $2$ categories.

First we define the following sheaf of categories $C_1$ on $M_1$:

To every $1-$connected open set $U$ of $M_1$ of affine chart, we
associate the category of affinely locally trivial affine bundles
over $U$ with typical fiber $(F_1,\nabla_{F_1})$ such that the
induced $\pi_{11}(U)$ bundle, is the restriction of
$\pi_{11}(M_1)$ to $U$.

We now define on $M_1$ the sheaf of $2$ categories $C_2$.

  For every open
set $U$, consider an element $e$ of $C_1(U)$. It is an affinely
locally trivial affine bundle over $U$ with typical fiber
$(F_1,\nabla_{F_1})$.
 We associate to $e$ the sheaf of category $C_2(e)$
which objects are a.l.t bundles over $e$ with typical fiber
$(F_2,\nabla_{F_2})$ such that on  $e$, the flat bundle with
typical fiber $Aff(F_2,\nabla_{F_2})/T_{F_2}$ induced, is the one
induced by the representation $\pi_1(F_1)\rightarrow
Aff(F_2,\nabla_{F_2})/T_{F_2}$. Gluing conditions for the $2$
sheaf $C_2$ are done using the representation
$\pi_1(M_1)\rightarrow Aut(\pi_{12}(F_1))$.

This enables to associate to the open set $U$,  $C_2(U):=\cup
C_2(e),e \in C_1(U)$.

\medskip

We will see that  $C_2$ is in fact a $2$ gerbe. An important fact
is that this $2$ gerbe is defined recursively, that is, we have
first define the gerbe $C_1$. This will considerably simplify the
expression of the classifying $3$ cocycle which in this case will
be an usual Cech $3$ cocycle.

\medskip

Let's precise now what $1$ and $2$ arrows are in the $2$ category
$C_2$.

Let  $e_1$ and $e_2$ be objects of $C_1(U)$ where $U$ is a simply
connected open set of affine chart. A map between $e_1$ and $e_2$
can be represented by an affine map $t:U\rightarrow T_{F_1}$
acting on $U\times (F_1,\nabla_{F_1})$ as follows:

$$
U\times (F_1,\nabla_{F_1})\longrightarrow U\times
(F_1,\nabla_{F_1})
$$
$$
(x,y)\longmapsto (x,t(x)y).
$$

The map $t$ can be viewed as an automorphism of the site of open
sets of $U\times F_1$, so the map $t$  induces a functor between
the sheaf of categories  $C_2(e_1)$ and $C_2(e_2)$, such maps $t$
are $1$ arrows of our $2$ category.

A $2$ arrow of our $2$ category $C_2(U)$ over the one arrow $t$
can be represented as a family of affine arrows
$$
t_{i_1}:U\times U_{i_1}\times (F_2,\nabla_{F_2})\longrightarrow
U\times t(U_{i_1})\times (F_2,\nabla_{F_2})
$$
$$
(x,y,z)\longmapsto (x,t(x)y,t_1(x,y)z),
$$
where $(U_{i_1})$ is an open covering of $(F_1,\nabla_{F_1})$ by
one connected open set of affine charts, and $t_1:U\times
U_{i_1}\rightarrow T_{F_2}$ is an affine map.

It is easy to see that our category satisfy the axioms which
defines $2$ gerbes, so we have:

\medskip

{\bf Proposition 1.2.} {\it The $2-$ sheaf $C_2$ is a $2$ gerbe.}

\medskip

{\bf The classifying three cocycle.}

\medskip

The representation $\pi_1(M_1)\rightarrow
Aff(F_1,\nabla_{F_1})/T_{F_1}$ induces an affine  flat bundle
$V_1$ over $(M_1,\nabla_{M_1})$ with typical fiber $T_{F_1}$, we
will call $S_1$ the sheaf of affine sections of this bundle. The
representation $\pi_1(F_1)\rightarrow
Aff(F_2,\nabla_{F_2})/T_{F_2}$ defines a
 flat bundle $V_{12}$ over $(F_1,\nabla_{F_1})$ with typical fiber $T_{F_2}$; we will
 call $S_{12}$ the sheaf of affine sections of this bundle.

  The representation
 $\pi_1(M_1)\rightarrow Aut(\pi_{12})(F_1)$ induces
 the sheaf $S_{123}$ of affine maps
 $U\rightarrow S_{12}$, where $U$ is an open set of $M_1$. This sheaf is a locally
 constant sheaf over $M_1$.

 \medskip

 Now we consider an open covering $(U_i)_{i\in I}$ of $(M_1,\nabla_{M_1})$
 by $1$ connected affine charts, For each $i$ the $2$ category
 $C_2(U_i)$ is not empty.
 We will choose in each $U_i$ an element $e_i$ of $C_1(U_i)$.
 Let $i,j$ such that $U_i\cap U_j$ is not empty. The restriction of
 $e_i$ and $e_j$ to $U_i\cap U_j$ gives rise to an arrow

 $\phi_{ij}: {e_i}_{\mid U_i\cap U_j}\rightarrow {e_j}_{\mid U_i\cap U_j}$.

 This arrow can be expressed as a map
 $$
 U_i\cap U_j\times (F_1,\nabla_{F_1})\longrightarrow
 U_i\cap U_j\times (F_1,\nabla_{F_1})
 $$
 $$
 (x,y)\longmapsto (x,t_{ij}(x)y).
 $$

Recall that the map $\phi_{ij}$ can also be viewed as a functor
$C_2(e_i)_{\mid U_i\cap U_j}\rightarrow C_2(e_j)_{\mid U_i\cap
U_j}$.

Now we consider the restriction of the functor $\phi_{ij}\circ
\phi_{jk}\circ \phi_{ik}^{-1}=\psi_{ijk}$ to  the family of
$U_i\cap U_j\cap U_k$.

It can be represented by a family of maps
$$
U_i\cap U_j\cap U_k\times U_{i_1}\times
(F_2,\nabla_{F_2})\longrightarrow U_i\cap U_j\cap U_k\times
t_{ij}t_{jk}t_{ik}^{-1}(U_{i_1})\times (F_2,\nabla_{F_2})
$$
$$
(x,y,z)\longmapsto (x, t_{ij}(x)
t_{jk}(x)t_{ki}(x)(y),u_{ijk}(x,y)(z)),
$$
where $U_{i_1}$ is a one connected open set of affine chart of
$(F_1,\nabla_{F_1})$, $u_{ijk}:U_i\cap U_j\cap U_k\times
U_{i_1}\rightarrow T_{F_2}$ is an affine map.

\medskip

Now, we can restrict $\psi_{ijk}$ to $U_i\cap U_j\cap U_k\cap
U_l=U_{ijkl}$, writing the boundary of the chain $\psi_{ijk}$, we
obtain $\rho_{ijkl}=
\psi_{jkl}\psi_{ikl}^{-1}\psi_{ijl}\psi_{ijk}^{-1}$, $\rho_{ijkl}$
can be viewed as a map
$$
U_{ijkl}\times U_{i_1}\times (F_2,\nabla_{F_2}) \longrightarrow
U_{ijkl}\times U_{i_1}\times (F_2,\nabla_{F_2})
$$
$$
(x,y,z)\longmapsto (x,y,w_{ijkl}(x,y)z),
$$
as the family of map $t_{ij}t_{jk}t_{ki}$ defines a $2$ Cech
cocycle of $S_1$. The family $\rho_{ijkl}$ can be viewed as
sections of the bundle $S_{123}$.

\medskip

{\bf Theorem 1.3.} {\it The family of maps $\rho_{ijkl}$ that we
have just define is a $3$ Cech cocycle.}

\medskip

{\bf Proof.}

 We must calculate the boundary of the family of $\rho_{ijkl}$.

 Let $U_{ijklm}$ be the intersection
 $U_i\cap U_j\cap U_k\cap U_l\cap U_m$, we have:

 $d(\rho_{ijklm})= $

 $\rho_{jklm}-\rho_{iklm}+\rho_{ijml}-\rho_{ijkm}+\rho_{ijkl}$

$=\psi_{klm}\psi_{jlm}^{-1}\psi_{jkm}\psi_{jkl}^{-1}$

$-(\psi_{klm}\psi_{ilm}^{-1}\psi_{ikm}\psi_{ikl}^{-1})$

$+\psi_{jml}\psi_{iml}^{-1}\psi_{ijl}\psi_{ijm}^{-1}$

$-(\psi_{jkm}\psi_{ikm}^{-1}\psi_{ijm}\psi_{ijk}^{-1})$

$+\psi_{jkl}\psi_{ikl}^{-1}\psi_{ijl}\psi_{ijk}^{-1}=0.$

\medskip

The associated $3$ cocycle $\rho_{ijkl}$ is not the obstruction to
the existence of a composition serie, suppose that it vanishes.

This means that there is a family of maps
$$
h_{ijk}:U_{ijk}\times U_{i_1}\times (F_2,\nabla_{F_2})
\longrightarrow U_{ijk}\times \psi_{ijk}(U_{i_1})\times
(F_2,\nabla_{F_2})
$$
$$
(x,y,z)\longmapsto (x,\psi_{ijk}(x)y,h'_{ijk}(x,y)z),
$$
(where the map $h'_{ijk}: U_{ijk}\times U_{i_1}\rightarrow
T_{F_2}$ is an affine map which boundary is $\rho_{ijkl}$) which
 is a $2$ cocycle of $S_{1}\oplus S_{123}$.

 We have:

\medskip

{\bf Theorem 1.4.} {\it If the cocycle $\rho_{ijkl}$ is trivial,
then two cocycle $h_{ijk}$ that we have just define is the
obstruction to the existence of a composition serie associated to
the bundles $S_1$, $S_{12}$ and $S_{123}$.}

\medskip

{\bf Proof.}

If the cocycle $h_{ijk}$ is trivial, then there exists a family of
maps $b_{ij}: U_i\cap U_j\longrightarrow S_{1}\oplus S_{123}$ such
that the family of map $(t_{ij}+b_{ij})$ define a $1$ Cech
cocycle. This implies that the family of maps $t_{ij}$ is a $1$
cocycle up to a $1$ boundary, then it defines an affine bundle
over $(M_1,\nabla_{M_1})$ with typical fiber $(F_1,\nabla_{F_1})$
associated to $\pi_{11}$. Let $(M_2,\nabla_{M_2})$ be its total
space. The  obstruction of the existence of an affine bundle over
$(M_2,\nabla_{M_2})$ with typical fiber
 $(F_2,\nabla_{F_2})$ associated to $S_{1}$, $S_{12}$ and $S_{123}$
 is  given by $h_{ijk}$, since the open set $U\times U_{i_1}$ used to build
 the obstruction $\rho$ can be viewed as open subsets of $M_2$.

 \medskip

 {\bf Remark.}

 When the obstruction $\rho$ vanishes, the cocycle $h_{ijk}$ define on
 $M_1$ a gerbe which can be viewed
 as trivial $2$ gerbe.

 \bigskip

 {\bf 2. The general case.}

\medskip

In this part, we will classify  composition series
$(M_n,\nabla_{M_n})\rightarrow...\rightarrow
(M_2,\nabla_{M_2})\rightarrow (M_1,\nabla_{M_1})$.

We will denote by $(F_i,\nabla_{F_i})$ the fiber of the a.l.t
affine bundle $f_i:(M_{i+1},\nabla_{M_{i+1}})\rightarrow
(M_i,\nabla_{M_i}).$

Let $i<j\leq n$, the map $f_{ij}=f_{j-1}\circ f_{j-2}\circ...\circ
f_i,$ is an affine map $f_{ij}:(M_j,\nabla_{M_j})\rightarrow
(M_i,\nabla_{M_i}).$ Since this map is a submersion and $M_j$ is
compact,
 we deduce that $(M_j,\nabla_{M_j})$ is
the total space of a locally trivial differentiable fibration over
$(M_i,\nabla_{M_i})$. The Serre bundle theorem implies the
following
 exact sequence

$$
1\rightarrow \pi_1(F_i)\rightarrow
\pi_1(M_{i+1})\rightarrow\pi_1(M_{i}) \rightarrow 1
$$
when $j=i+1$.

The affine map $f_{ij}$ define an affine bundle which is not
necessarily an affinely locally trivial affine bundle.

\medskip

Recall that the map $f_i$ gives rise to a representation
$\pi_{i}:\pi_1(M_i)\rightarrow Aff(F_i,\nabla_{F_i})/T_{F_i}$, and
to a flat bundle $V_{i}$ over $M_i$ with typical $T_{F_i}$.

Thus we have a flat bundle $V_{n-1}$ over $M_{n-1}$ with typical
fiber $T_{F_{n-1}}$. The group $\pi_1(F_{n-2})$ is a subgroup of
$\pi_1(M_{n-1})$, thus we have a flat bundle $V_{n-1n-2}$ over
$F_{n-2}$ with typical fiber $T_{F_{n-1}}$ induced by $V_{n-1}$.

The fundamental group $\pi_1(M_{n-2})$ of $M_{n-2}$, acts on
$V_{n-1n-2}$ via the representation $\pi_{n-1}$. This action
defines a flat $V'_{n-1n-2}$ over $M_{n-2}$ with typical fiber
$V_{n-1n-2}$. This bundle also gives rise to a bundle
$V_{n-1n-2n-3}$ over $F_{n-3}$ as $\pi_1(F_{n-3})$ is a subgroup
of $\pi_1(M_{n-2})$. Recursively, we can define bundle
$V_{n-1...n-i}$ over $n-i$ with typical fiber $V_{n-1...n-i+1}$
over $F_i$.

We can also define the representation $\pi_{n-1n-2}$ which is the
restriction of $\pi_{n-1}$ to $\pi_1(F_{n-2})$. This
representation induces on $F_{n-2}$ a flat bundle $\pi_{n-1n-2}'$
with typical fiber $Aff(F_{n-1},\nabla_{F_{n-1}})/T_{F_{n-1}}$,
$\pi_1(M_{n-2})$ acts on this bundle via $\pi_{n-1}$, one deduces
a flat bundle over $M_{n-2}$ with typical fiber $\pi_{n-1n-2}'$
which induce a flat bundle $\pi_{n-1n-2n-3}'$ over $F_{n-3}$.
Recursively, we can define bundle $\pi_{n-1...1}'$.

Remark also that considering the composition serie
$(M_j,\nabla_{M_j})\rightarrow....\rightarrow (M_1,\nabla_{M_1})$
for $j<n$, one can also define the bundle $V_{jj-1....i}$ with
typical fiber $V_{j...i+1}$ and base space $F_i$.

\medskip

Let Let $S_{n-1n-2}$ be the sheaf of affine sections of
$V_{n-1n-2}$ one may define the sheaf $S_{n-1n-2n-3}$ of affine
sections of $S_{n-1n-2}$ over $F_{n-3}$. Recursively, we can also
define the sheaf $S_{n-1...1}$ of affine sections of $S_{n-1...2}$
over $M_1$. The gluing conditions for those sheaves are given by
the bundles $\pi_{n-1...i}'$. One can also define in the same way
the bundles $S_{i...1}$, $i\leq 1$.

 \medskip

 {\bf Remark.}

 The bundles $V_{j...1}$ and $S_{j...1}$ depend only of the affine structures of
 $(F_1,\nabla_{F_1})$,...,$(F_j,\nabla_{F_j})$ and $(M_1,\nabla_{M_1})$.
 They can be defined without suppose the existence of a
 composition serie.

 \medskip

 {\bf The classifying $n-$cocycle.}

 \medskip

 Given bundles $V_{ii-1...1}$, and $S_{i...1}$ $1\leq i\leq n-1$ as above,

  We want to classify all
 composition series $(M_n,\nabla_{M_n})\rightarrow...\rightarrow (M_1,\nabla_{M_1})$
 associated
 to the family of bundles
 $V_{i...1}$. To make this classification, we are first going to define an
 $n-$cocycle.

 \medskip

 First we consider the trivialization of the flat bundle over $M_1$ with typical
 fiber $Aff(F_1,\nabla_{F_1})/T_{F_1}$ over $M_1$ induced by $\pi_1$.

 It is defined by
 $$
 U_i\cap U_j\times Aff(F_1,\nabla_{F_1})/T_{F_1}
 \longrightarrow U_i\cap U_j\times Aff(F_1,\nabla_{F_1})/T_{F_1}
 $$
 $$
 (x,y)\longmapsto (x,\bar t_{ij}(y))
 $$
 Consider for each $x$ in $U_i\cap U_j$
 an element  $t_{ij}(x)$ of $Aff(F_1,\nabla_{F_1})$
 over $\bar t_{ij}$ which depends affinely of $x$,
  one may define the $2$ cocycle
  $$
  t_{ijk}:U_i\cap U_j\cap U_k\longrightarrow T_{F_1}
  $$
  $$
  x\longmapsto t_{ij}(x)t_{jk}(x)t_{ki}(x).
  $$
  The family of $t_{ijk}$ may be considered as local sections of
  the bundle $S_1$. It is the cocycle associated to the gerbe
 which at $U_i$ associated the category of a.l.t affine bundles
 over $U_i$ with typical fiber $(F_1,\nabla_{F_1})$, such that
 the bundle over $U_i$ with typical fiber $T_{F_1}$
 associated is the restriction of $V_1$.

  The map $t_{ijk}$ induces a functor on the category of open sets of
  $U_i\cap U_j\cap U_k\times (F_1,\nabla_{F_1})$.

  Consider a covering
  $U_{i_1}$ of $(F_1,\nabla_{F_1})$ by affine $1$ connected affine charts.

  Let $U_{ijk}=U_i\cap U_j\cap U_k$.
  We  associate to $U_{ijk}\times U_1$ the category of a.l.t affine bundles
  with typical fiber $(F_2,\nabla_{F_2})$ such that the bundle with typical
  fiber $T_{F_2}$ associated is induced $V_2$.

  On $U_{ijkl}=U_i\cap U_j\cap U_k\cap U_l$,  the boundary of
  $t_{ij}t_{jk}t_{ki}$
  gives rise to the map
  $$
 v_{ijkl}: U_{ijkl}\times U_{i_1}\times (F_2,\nabla_{F_2})\longrightarrow
  U_{ijkl}\times U_{i_1}\times (F_2,\nabla_{F_2})
  $$
  $$
  (x,x_1,x_2)\longmapsto (x,x_1,u_{ijkl}(x,x_1)(x_2)),
  $$
where the map $u_{ijkl}$ are affine sections of the bundle
$S_{321}$.
 The family of maps $u_{ijkl}$ is a $3$ cocycle.

\medskip

Let $U_{1...i}=U_1\cap U_2....\cap U_i$. For each $j$, we will
consider an open covering $U_{i_j}$ of $(F_j,\nabla_{F_j})$ by
$1-$ connected open sets of affine charts.

Suppose that we have defined a family of maps
$$
v_{1...j}:U_{1...j}\times U_{i_1}\times...\times U_{i_{j-3}}\times
(F_{j-2},\nabla_{F_{j-2}}) \longrightarrow U_{1...j}\times
U_{i_1}\times...\times U_{i_{j-3}}\times
(F_{j-2},\nabla_{F_{j-2}})
$$
$$
(x,x_1,...,x_{j-2})\longmapsto
(x,x_1,...,x_{j-3},u_{1...j-3}(x,x_1,...x_{j-3}) (x_{j-2}))
$$
(where the family of maps $u_{1...j-3}$ are local affine sections
of the bundle $S_{j...1}$) which represent a $j-1$ Cech cocycle.
Then on $U_{1...j}\times U_{i_1}\times...\times
(F_{j-2},\nabla_{F_{j-2}})$, one can define the sheaf of category
$C_j$ such that the objects $C_j(U_{1...j}\times...\times
U_{i_{j-2}})$ are affinely locally trivial affine bundles with
typical fiber $(F_{j-1},\nabla_{F_{j-1}})$ over $U_{1...j}\times
U_{i_1}\times...\times (F_{j-2},\nabla_{F_{j-2}})$, and the
canonical flat vector bundle with typical fiber $T_{F_{j-1}}$
associated is the restriction of $V_{j-1}$.

The map $v_{1...j}$ induces a functor $w_{1...j}$ in the category
$C_j(U_{i...1}\times...\times U_{i_{j-2}}).$

The restriction of the composition $w_{2...j+1}\circ...\circ
w_{1...\hat k...j+1}^{(-1)^{k+1}}\circ...\circ
w_{1...j}^{(-1)^{j+2}}$ on $U_{1...j+1}\times
U_{i_1}\times...U{i_{j-2}}$ induces a map
$$
v_{1...j+1}:U_{1...j+1}\times U_{i_1}...\times U_{i_{j-2}}\times
(F_{j-1},\nabla_{F_{j-1}}) \longrightarrow U_{1...j+1}\times
U_{i_1}...\times U_{i_{j-2}} \times (F_{j-1},\nabla_{F_{j-1}})
$$

$$
(x,x_1,...,x_{j-2})\longmapsto (x,x_1,...,x_{j-2},u_{1...{j-2}}
(x,...,x_{j-2})(x_{j-1}))
$$

\medskip

{\bf Proposition 2.1.} {The Cech chain $v_{1...j+1}$ that we have
just define recursively is a $j$ Cech cocycle.}

\medskip

{\bf Proof.}

The proof will be made recursively. We have already verify the
result if $j=2$.

Suppose that the chain $v_{1...j}$ is a Cech cocycle for $k\leq
j$, then the writing the boundary $d(v_{1...j+1})$ of
$v_{1...j+1}$, we obtain
$$
\sum_{l=1}^{l=j+2} (-1)^{l+1}v_{1...\hat l...j+2}
$$
$=$
$$
\sum_{l=1}^{j+2} (-1)^{l+1}(\sum_{m=1}^{m=l-1}
 w_{1..\hat m..\hat l..j+2}^{(-1)^{m+1}}
 \circ...\circ w_{1..\hat {l-1}\hat l...j+1}^{(-1)^{l}}
 +\sum_{m=l+1}^{j+2}w_{1..\hat l..\hat m..j+2}^{(-1)^{m}})=0.
$$

\medskip

The cocycle $v_{1...n+1}$ is not the obstruction of the existence
of a composition sequence associated to the family of bundles
$S_{ii-1...1}$. If its cohomology class is zero, its means that
there exists a chain $a_{1...n}$ which boundary is $v_{1...n+1}$.
Suppose that the chain $z_{n-1}=a_{1...n}+v_{1...n}$ considered as
an element of the Whitney sum of the bundles $S_{n-1...1}\oplus
S_{n-2...1}=T_{n-1n-2}$ is an $n-1$ cocycle.

If the cohomology class of the cocycle $z_{n-1}$ is zero, then it
is the boundary of an $n-2$ cocycle $a_{1...n-1}$. We can define
the chain $z_{n-2}=a_{1...n-1}+v_{1...n-1}$ considered as an
element of $T_{n-1n-2n-3}=T_{n-1n-2}\oplus S_{n-3...1}$.

Suppose that we have define the cocycle $z_{n-i}$ and its
cohomology class is zero. It is the boundary of a chain
$a_{1...n-i}$.
 This means that the chain $z_{n-(i+1)}=a_{1...n-i}+v_{1...n-i}$
is an $n-(i+1)$ cocycle viewed as an element of
$T_{n-1...n-(i+2)}=T_{n-1...n-(i+1)}\oplus S_{n-(i+2)...1}$.

\medskip

{\bf Theorem 2.2.} {\it Suppose that we can construct a chain
$z_2$ by the processus that we have just describe; then it is the
obstruction to the existence of a composition serie associated to
the family $S_{ii-1...1}$. When it vanishes the set of equivalence
classes of composition series is given by $H^1(M_1,T_{n-1...1})$.}

\medskip

{\bf Proof.}

Suppose that we can define the cocycle $z_2$ and its cohomology
class is zero, then then it is the boundary of an element $z_1$ of
$T_{n-1...1}$.

This implies that up to a boundary the chain $t_{ij}t_{jk}t_{ki}$
is zero,
 we deduce that we can
solve the first extension problem, there exists a bundle
$(M_2,\nabla_{M_2})\rightarrow (M_1,\nabla_{M_1})$ associated to
$S_1$.

The obstruction of the existence of a composition serie
$(M_3,\nabla_{M_3})\rightarrow (M_2,\nabla_{M_2})\rightarrow
(M_1,\nabla_{M_1})$ is also given by the class of $z_2$ which is
zero, so we can solve the second extension problem.

Suppose that we can solve the extension problem
$(M_i,\nabla_{M_i})\rightarrow...\rightarrow (M_1,\nabla_{M_1})$
then obstruction to solve the extension problem
$(M_{i+1},\nabla_{M_{i+1}})\rightarrow...\rightarrow
(M_1,\nabla_{M_1})$ is also given by the cohomology class of $z_2$
which is zero.

Recursively, we deduce that we can solve the extension problem
$(M_n,\nabla_{M_{n}})\rightarrow...\rightarrow
(M_1,\nabla_{M_1})$.

\bigskip

\bigskip

{\bf 3. The Conceptualization.}

\medskip

The resolution of the extension problem, given the base space
$(M_1,\nabla_{M_1})$ and the bundle $S_1$ has been done using the
gerbe theory. It is natural to think that the existence of a
composition serie $(M_n,\nabla_{M_n})\rightarrow ...\rightarrow
(M_1,\nabla_{M_1})$ must be solve using $n-1-$gerbe theory. In
this part, we are going to build a commutative $n-$gerbe theory.

\medskip

\medskip

{\bf On $n-$ Categories.}

\medskip

Supposed that the notion of  $i$ category is defined. An $i+1$
category  $C_{i+1}$, is given by

The class of objects $O(C_{i+1})$,

$i_1$

the morphisms $Hom_{C_{i+1}}(x,y)$ between two objects $x$ and $y$
of $C_{i+1}$ which is an $i$ category,

$i_2$

 For objects $x,y,z$ in $C_{i+1}$, the composition $i$ functor

$o_i:Hom_{C_{i+1}}(x,y)\times Hom_{C_{i+1}}(y,z) \rightarrow
Hom_{C_{i+1}}(x,z)$

$i_3$

 For each objects $x_1,...,x_{i+4}$  in $C_{i+1}$, we will assume that
the following strict condition: the  composition

$$
Hom_{C_{i+1}}(x_{i+3},x_{i+4})\times...\times
Hom_{C_{i+1}}(x_2,x_3)\times Hom_{C_{i+1}}(x_1,x_2)\rightarrow
Hom_{C_{i+1}}(x_1,x_{i+4})
$$
does not depend of the order in which it is made.
 More others conditions need to be specified, but we don't need
them.

\medskip

{\bf The notion of sheaf of $i$ categories.}

\medskip

Now we define recursively the notion of sheaf of $i$ categories on
a topological space $M$.

We assume known the notion of sheaves of sets.

Supposed that we have already defined the notion of  sheaves of
$i$ categories. An  sheaf of $i+1$ categories $C_{i+1}$ on the
topological space $M$, will be a map which assign to each open set
$U$ a $i+1$ category $C_{i+1}(U)$, such that

For each inclusion $U\hookrightarrow V$, there exists a $i+1$
functor

$$
c_{U,V}:C_{i+1}(V)\rightarrow C_{i+1}(U)
$$
which satisfies $c_{U,V}\circ c_{V,W}=c_{U,W}$ for any open sets
$U,V$ and $W$ such that $U\hookrightarrow V\hookrightarrow W$.

Gluing condition for objects. Consider an open covering $(U_j)$ of
$M$, and in each $U_i$ an object $A_j$, we restrict each object We
restrict each object of the family $(A_{j_1},...,A_{j_{i+3}})$ to
$U_{j_1...j_{i+3}}$, we assume that if the composition
$$
Hom_{C(U_{j_1...j_{i+3}})}(A_{j_{i+2}},A_{j_{i+3}})\times...\times
Hom_{C(U_{j_1...j_{i+3}})}(A_{j_{1}},A_{j_{2}})\longrightarrow
Hom_{C(U_{j_1...j_{i+3}})}(A_{j_1},A_{j_{i+3}})
$$
does not depend on the order in which it is made then there exists
a global object $A$ which restriction on each $U_i$ is $A_i$.

Gluing condition for arrows. For each pair of global objects $x$
and $y$,  $Hom_{C_{i+1}(M)}(x,y)$ is a sheaf of $i$ categories.

\medskip

{\bf The notion of $n$ gerbe.}

\medskip

Now, we will define the notion of $n$ gerbe where $n\geq 2$.

Consider a topological space $M$,  endowed with a sheaf of $i$
categories $C_i$ for each $i=1,...,n$ such that

$g_1$

For each open set $U$, the set of objects of the category $C_i(U)$
is the same for each $i=1,...n$.

$g_2$

For each $x\in M$, we suppose that there is a neighborhood $U_x$
of $x$ such that $C_i(U_x)$ is not empty.

$g_3$

 The sheaf of categories category $C_1$ is a gerbe with lien the abelian sheaf $T_1$ over
$M$.

$g_4$

 The    arrows between two objects $x, y$ considered as elements of
the $2$ category $C_2(U)$ is a category $Hom_{C_2(U)}(x,y)$ which
objects are elements of $Hom_{C_1(U)}(x,y)$ the  arrows between
the objects $x$ and $y$ considered as elements of $C_1(U)$. The
composition $o_2:Hom_{C_2(U)}(x,y)\times
Hom_{C_2(U)}(y,z)\rightarrow Hom_{C_2(U)}(x,z)$ transforms two $1$
arrows $f$ and $g$ to $gf$ where the product $gf$ is considered in
respect to the one of $Hom_{C_1(U)}(x,y)\times
Hom_{C_1(U)}(y,z)\rightarrow Hom_{C_1(U)}(x,z)$.

We recall that for each $x$ and $y$ in $C_2(U)$ the $1$ arrows
between $x$ and $y$ are the objets of the category
$Hom_{C_2(U)}(x,y)$, we have just precise how the functor $o_2$
acts on objects, not how it acts on maps. We deduce that the
product is known up to $2$ arrows.

The $1$ arrows between $x$ and $y$ considered as elements of
$C_3(U)$ is $Hom_{C_1(U)}(x,y)$, the $2$ arrows are the $2$ arrows
of $Hom_{C_2(U)}(x,y)$.

The product $o_3:Hom_{C_3(U)}(x,y)\times
Hom_{C_3(U)}(y,z)\rightarrow Hom_{C_3(U)}(x,z)$ transform  two $1$
arrows $f$ and $g$ onto $gf$, the product is considered in respect
to the one of $C_1(U)$ and two $2$ arrows $h$ and $k$ in $kh$ the
product, considered is the one of $C_2(U)$. We can also say that
the product is known up to $3$ arrows.

Recursively, suppose that we have defined the $1,...,i$ arrows of
the category $C_i(U)$, then the   $1$ arrows of $C_{i+1}(U)$ are
the arrows of $C_1(U)$,..., the  $i$ arrows of the category
$C_{i+1}(U)$ is the $i$ arrows of the category $C_i(U)$.

$g_5$

Suppose also that we have define recursively the product of $l\leq
i$ arrows of $Hom_{C_i(U)}(x,y)$. Then the product $o_{i+1}:
Hom_{C_{i+1}(U)}(x,y)\times Hom_{C_{i+1}(U)}(y,z)\rightarrow
Hom_{C_{i+1}(U)}(x,z)$ is an $i$ functor which send
 two $l$ arrows $l\leq i$ in $C_{i+1}(U)$ onto the
one  with  respect to the composition in $C_l(U)$. We can also
remark that the product $o_{i+1}$ is defined up to $i+1$ arrows.

$g_6$

We suppose that in $C_1(U)$ the arrows are invertible, in
$C_2(U)$, a $1$ arrow is invertible up to a $2$ arrow, a $2$ arrow
is invertible, in $C_i(U)$, a $1$ arrow is invertible up to a $2$
arrow, a $2$ arrow is invertible up to a $3$ arrow,... a $i-1$
arrow is invertible up to a $i$ arrow, and  $i$ arrow are
invertible.

$g_7$

 Given an object $f$ of $Hom_{C_2}(x,x)$,
where $x$ is an object of $C_2(U)$, the set of morphisms of $f$ is
isomorphic to  $T_2(U)$ where $T_2$ is a sheaf over $M$. More
generally,  let $g$ be an $i-1$ map in the category $C_i(U)$, the
set of morphisms of $g$ is isomorphic to $T_i(U)$ where $T_i$ is a
sheaf over $M$.

$g_8$

Any two objets of $C_i(U)$ $i\leq n$ are locally isomorphic.

\medskip

{\bf Definition 3.1.}

A family of sheaves of $i$ categories $C_i$ for each $i\leq n$
which satisfy the conditions $g_1,...,g_8$, will be called an
$n-$gerbe.

\medskip

\bigskip

{\bf The classifying $n+1-$cocycle.}

\medskip

In this part, we are going to consider $n$ gerbes such that for
each $i$, the sheaf $T_i$ is a commutative sheaf.

 We are going to associate to each $n-$gerbe a classifying $n+1-$Cech cocycle
 which takes values in $T_n$.

We will assume that $M$ is a manifold, $U_i$ is an open covering
of $M$ such that $C_i(U)$ is not empty and for each family
$\{i_1,...,i_k\}$ we set $U_{i_1...i_k}=U_{i_1}\cap U_{i_2}...\cap
U_{i_k}$.

\medskip

We first choose in each open subset $U_i$ an object $x_i$, we may
consider $x_i$ as an object of $C_1(U_i)$. If $U_{ij}$ is not
empty, then there is an arrow of $t^1_{ij}$ between $x_i$ and
$x_j$ (which is an isomorphism). we can write the two chains
$$
c_2(U_{i_1i_2i_3})=t^1_{i_1i_2}t^1_{i_2i_3}t^1_{i_3i_1}.
$$
As $C_1$ is a gerbe, this implies that the Cech boundary of $c_2$
is trivial.

But one may also consider arrows of $C_1$ as $1$ arrows of $C_2$,
this implies that the $d(c_2)$ is trivial up two a $3-$ $T_2$
chain $c_3(U_{i_1i_2i_3i_4})$ represented by a $3$ arrows. It
results from $i_3$ that $c_3$ is a $3-$ Cech cocycle.

Suppose that we have defined recursively a  Cech $j+1-$cocycle
associated to the sheaf of $j-$categories $C_j$. It is a family of
elements $c_{j+1}(U_{i_1...i_{j+2}})$ of $j$ arrows of
$C_j(U_{i_1...i_{j+2}})$. As the boundary of $c_{j+1}$ is trivial,
it may be represented by an $j+2$  Cech chain $c_{j+2}$ of
$T_{j+1}$ which is a $j+2$ cocycle.
 So we can deduce recursively the existence of an $n+1$ $T_n$ cocycle $c_{n+1}$
 associated
 to the $n-$ gerbe.

\medskip

{\bf Definitions 3.2.}

An $n-$gerbe, will be said $n-$trivial if the cohomology class of
the associated $n+1-$cocycle $c_{n+1}$ that we have just define is
trivial.

\medskip

Let $C$, and $C'$ two $n-$gerbes associated to the family of
sheaves $T_1$,..., $T_n$, which are locally isomorphic i.e,  each
$x$ in $M$, has an open neighborhood $U_x$ such that $C_i(U_x)$
and $C_i'(U_x)$ are not empty, and
 objects of $C_i(U_x)$ and $C'_i(U_x)$ are locally
isomorphic. We will say that the locally isomorphic $n-$gerbes $C$
and $C'$ are  equivalent if and only if for every open set $U$ of
$M$, there is a $n$ isomorphism $\phi_n(U):C_n(U)\rightarrow
C'_n(U)$ such that the restrictions of $\phi_n(U)$ and $\phi_n(V)$
on $U\cap V$ coincide with $\phi_n(U\cap V)$.

Suppose that the class  $[c_{n+1}]$ of the cocycle $c_{n+1}$ is
trivial. It implies that it is the boundary of a chain $a_{n}$.
Let consider the chain $z_n=a_{n}+c_{n}$ of $T_{n-1}\oplus T_n$,
if its  cohomology class is trivial, it implies that it is the
boundary of a chain $a_{n-1}$, we set $z_{n-1}=a_{n-1}+c_{n-1}$.
Suppose that recursively we have defined the chain $z_{n-i}$.

\medskip

 The $n$ gerbe is said $n-i$ trivial, if we can define the cocycle
 $z_{n-i}$ by the processus above. The fact that the $n$ gerbe is
 $n-i$ trivial means that it can be considered as a $n-i-1$ gerbe.
 as follows: we consider the $n-i-1$ gerbe $C'$ such that for every
 open set $U$, $C'_j(U)=C_j(U)$
if $j<n-i-1$. We define the $n-i-1$ maps of $C_{n-i-1}(U)$ to be
$T_{n-i-1}(U)\oplus...\oplus T_n(U)$.

\medskip

{\bf Proposition 3.3.} {\it

The set of  equivalence classes of locally isomorphic  $n$ trivial
gerbes is given by $H^{n+1}(M,T_n)$ .}

\medskip

{\bf Proof.}

We have assigns to every $n-$gerbe a $n+1$ cocycle $c_{n+1}$ if it
class is trivial, it is equivalent to saying that the $n-$gerbe is
trivial. This implies that the map $C\rightarrow c_{n+1}$ is
injective.

On the other hand, given a cocycle $c_{n+1}$, we can find a family
of cocycles $c_2$,...,$c_n$, such that $c_{i+1}$ is induced by
$c_i$ by the processus described to build the classifying cocycle,
and $c_n$ induces $c_{n+1}$.

\bigskip

We can extend our definition and defines $\infty-$gerbe.

\medskip

{\bf Definition 3.4.}

 The family $C_n$ $n\in{\N}$
of $n-$sheaves of categories over the manifold $M$ is an
$\infty-$gerbe $C$ if and only if for each $n\in {\N}$, the family
$C_1$,...,$C_n$ is an $n-$gerbe and the family of sheaves $T_n$,
$n\in {\N}$ is an inductive system such that the map
$i_n:T_n\rightarrow T_{n+1}$ sends $c_{n+1}$ onto $c_{n+2}$, where
$c_{n+1}$ is the classifying cocycle associated to the gerbe
$C_1$,... $C_n$. We will call the inductive limit of $c_n$ the
classifying cocyle.

\bigskip
\bigskip

{\bf Acknowledgements.}

\medskip

The author would acknowledge the Abdus Salam ICTP, Trieste, Italy
for support.  This work has been done
 at ICTP.

\bigskip

\centerline{\bf Bibliography.}

\bigskip

[Bo] Borel, A. Linear algebraic groups. W.A. Benjamin, Inc, New
York-Amsterdam 1969.

[Bre] Breen, L. On the classification of $2-$gerbes and
$2-$stacks. Asterisque, 225 1994.

[Br] Bredon, G. E. Sheaf theory. McGraw-HillBook Co., 1967.

[Bry] Brylinski, J.L Loops spaces, Characteristic Classes and
Geometric Quantization, Progr. Math. 107, Birkhauser, 1993.

[Br-Mc] Brylinski, J.L, Mc Laughlin D.A, The geometry of degree
four characteristic classes and of line bundles on loop spaces I.
Duke Math. Journal. 75 (1994) 603-637.

[Ca] Carriere, Y. Autour de la conjecture de L. Markus sur les
vari\'et\'es affines. Invent. Math. 95 (1989) 615-628.

[De] Deligne, P. Theorie de Hodge III, Inst. Hautes Etudes Sci.
Publ. Math. 44 (1974), $5-77$.

[Du] Duskin, J. An outline of a theory of higher dimensional
descent, Bull. Soc. Math. Bel. S\'erie A 41 (1989) 249-277.

[F] Fried, D. Closed similarity affine manifolds. Comment. Math.
Helv. 55 (1980) 576-582.

[F-G] Fried, D. Goldman, W. Three-dimensional affine
crystallographic groups. Advances in Math. 47 (1983), 1-49.

[F-G-H] Fried, D. Goldman, W. Hirsch, M. Affine manifolds with
nilpotent holonomy. Comment. Math. Helv. 56 (1981) 487-523.

[G1] Goldman, W. Two examples of affine manifolds. Pacific J.
Math. 94 (1981) 327-330.

[G2] Goldman, W. The symplectic nature of fundamental groups of
surfaces. Advances in in Math. 54 (1984) 200-225.

[G3] Goldman, W. Geometric structure on manifolds and varieties of
representations. 169-198, Contemp. Math., 74

[G-H1] Goldman, W. Hirsch, M. The radiance obstruction and
parallel forms on affine manifolds. Trans. Amer. Math. Soc. 286
(1984), 629-949.

[G-H2] Goldman, W. Hirsch, M. Affine manifolds and orbits of
algebraic groups. Trans. Amer. Math. Soc. 295 (1986), 175-198.

[Gr] Grothendieck, A. Pursuing stacks, preprint
 avaible at University of Bangor.

God] Godbillon, C. Feuilletages. Etudes g\'eom\'etriques. Progress
in Mathematics, 98.

[K] Koszul, J-L. Vari\'et\'es localement plates et convexit\'e.
Osaka J. Math. (1965), 285-290.

[K] Koszul, J-L. D\'eformation des connexions localement plates.
Ann. Inst. Fourier 18 (1968), 103-114.

[Mc] Maclane, S. Homology. Springer-Verlag, 1963.

[Ma] Margulis, G. Complete affine locally flat manifolds with a
free fundamental group. J. Soviet. Math. 134 (1987), 129-134.

[Mi] Milnor, J. W. On fundamental groups of complete affinely flat
manifolds, Advances in Math. 25 (1977) 178-187.

[S-T] Sullivan, D. Thurston, W. Manifolds with canonical
coordinate charts: some examples. Enseign. Math 29 (1983), 15-25.

[T1] Tsemo, A. Th\`ese, Universit\'e de Montpellier II. 1999.

[T2] Tsemo, A. Automorphismes polynomiaux des vari\'et\'es
affines. C.R. Acad. Sci. Paris S\'erie I Math 329 (1999) 997-1002.

[T3] Tsemo, A. D\'ecomposition des vari\'et\'es affines. Bull.
Sci. Math. 125 (2001) 71-83.

[T4] Tsemo, A. Dynamique des vari\'et\'es affines. J. London Math.
Soc. 63 (2001) 469-487.

[T5] Tsemo, A. Fibr\'es affines to be published in Michigan J.
Math. vol 49.

\end{document}